\documentclass{article}
\usepackage{color}
\usepackage{amsmath}
\usepackage{amssymb}

\begin{document}

\title{Well posedness of an angiogenesis related integrodifferential diffusion model}
\author{Ana Carpio\thanks{Departamento de Matematica Aplicada, Universidad  Complutense, 28040 Madrid, Spain, 
ana\_carpio@ucm.es},
Gema Duro\footnote{Departamento de 
An\'alisis Econ\'omico: Econom\'{\i}a
Cuantitativa, Universidad Aut\'onoma de Madrid, 28049 Madrid, Spain}}

\date{May 2, 2015}

\maketitle

{\bf Abstract.} 
We prove existence and uniqueness of nonnegative solutions
for a nonlocal in time integrodifferential diffusion system related 
to angiogenesis descriptions. Fundamental solutions of 
appropriately chosen parabolic operators with bounded coefficients  
allow us to generate sequences of approximate solutions. 
Comparison principles and integral equations provide 
uniform bounds ensuring some convergence properties
for iterative schemes and providing stability bounds. 
Uniqueness follows from chained integral inequalities. \\


\section{Introduction}
\label{sec:intro}

Models for angiogenesis (blood vessel generation) have a complex mathematical structure, involving integral terms, transport operators, degenerate diffusion and eventually measure valued coefficients. In \cite{capasso}, for instance, the following model for the dynamics of the density of blood vessels is used to describe tumor induced angiogenesis:
\begin{eqnarray} 
\frac{\partial}{\partial t} p(t,\mathbf{x},\mathbf{v}) &=& 
 \alpha(c(t,\mathbf{x})) \rho(\mathbf{v}) p(t,\mathbf{x},\mathbf{v}) - \gamma p(t,\mathbf{x},\mathbf{v}) \int_0^t ds \, \tilde{p}(s,\mathbf{x})  \nonumber\\ 
&& - \mathbf{v}\cdot \nabla_\mathbf{x}   p(t,\mathbf{x},\mathbf{v}) + k \nabla_\mathbf{v} \cdot (\mathbf{v} p(t,\mathbf{x},\mathbf{v})) 
\nonumber\\ 
&& - \nabla_\mathbf{v} \cdot \left[\mathbf{F}\left(c(t,\mathbf{x})\right)
p(t,\mathbf{x},\mathbf{v})  \right]\!+  \sigma  
\Delta_\mathbf{v} p(t,\mathbf{x},\mathbf{v}), \label{vfpp} \\
\alpha(c)&=&\alpha_1\frac{\frac{c}{c_R}}{1+\frac{c}{c_R}} \geq 0, \quad
 {\bf F}(c)= \frac{d_1}{(1+\gamma_1c)^{q_1}}\nabla c,
\label{Falpha} \\
\frac{\partial}{\partial t}c(t,\mathbf{x})&=&d \Delta_{\mathbf x} c(t,\mathbf{x}) - \eta c(t,\mathbf{x})|\mathbf{j}(t,\mathbf{x})| \label{TAF}, \\
\mathbf{j}(t,\mathbf{x})&=& \int_{\mathbb R^2} \mathbf{v}' p(t,\mathbf{x},\mathbf{v}')\, d{\bf v}',
\quad \tilde{p}(t,\mathbf{x})= \int_{\mathbb R^2}  d{\bf v}' p(t,\mathbf{x},\mathbf{v}'),
\label{jp}
\end{eqnarray}
when $(\mathbf{x},\mathbf{v}) \in \mathbb{R}^2 \times \mathbb{R}^2$ and
$t \in[0, \infty)$.
{ This type of models is inspired in previous work on 
self-organized phenomena and retinal angiogenesis \cite{capasso2,capasso3}. }
Here, $p$ represents the density { (in space and
velocity)} of blood vessels growing in response to the distribution 
of the tumor angiogenic factor $c$. { Parameters}
$\sigma$, $\gamma$, $k$, $\alpha_1$, $c_R$, $d_1$, 
$\gamma_1$, $d$ and $\eta$ are positive constants. 
$\rho(\mathbf{v})$ is a gaussian profile centered at a point 
$\mathbf{v}_0$. 
{  The original model contains a Dirac mass 
$\delta_{\mathbf v_0}$ instead of a gaussian. Regularizations
of the form $\rho_{\varepsilon}(\mathbf v)= {1\over (\pi \varepsilon)^{N/2}} e^{-|\mathbf v- \mathbf v_0|^2 \over \varepsilon}$ 
are used for numerical purposes. The functions $\rho_{\varepsilon}$
tend to $\delta_{\mathbf v_0}$ as $\varepsilon$ tends to $0$.
$\rho(\mathbf v)$ stands for any of them.
From the theoretical point of view, we may seek to construct
solutions for delta valued coefficients as limits of solutions
for these gaussians approximations. However, proving compactness 
of such sequences of solutions, even in simpler models 
like the ones we consider here, is yet an open 
problem, as we will discuss in Section \ref{sec:discussion}.}

Establishing well posedness of these regularized problems is a challenging 
task due to the combination of degenerate diffusion, nonlinear transport terms and integrodifferential sources.  Here, we will use a simpler diffusion model as 
a basis to develop strategies to handle some of the technical difficulties:
\begin{eqnarray} \frac{\partial}{\partial t} p(t,\mathbf{x},\mathbf{v}) - \sigma
\Delta_{\mathbf{x}\mathbf{v}} p(t,\mathbf{x},\mathbf{v})  &=&  
 \alpha(c(t,\mathbf{x})) \rho(\mathbf{v})
 p(t,\mathbf{x},\mathbf{v}) \nonumber \\
& &- \gamma p(t,\mathbf{x},\mathbf{v}) \int_0^t \! \! d\,s \!\! 
\int_{\mathbb R^2}  \! d{\bf v}'  p(s,\mathbf{x},\mathbf{v}'), \label{heat1} \\
p(0,\mathbf{x},\mathbf{v})&=&p_0(\mathbf{x},\mathbf{v}), \label{heat2} \\
\frac{\partial}{\partial t}c(t,\mathbf{x}) - d \Delta_{\mathbf{x}} c(t,\mathbf{x}) 
&=& - \eta c(t,\mathbf{x}) { j(t,\mathbf{x})}, \label{heat3} \\
c(0,\mathbf{x})&=&c_0(\mathbf{x}), \label{heat4}
\end{eqnarray}
set in the whole space $(\mathbf{x},\mathbf{v}) \in \mathbb{R}^2 \times \mathbb{R}^2$, for  $t \in[0, \infty)$. { Here, 
\begin{eqnarray}
j(t,\mathbf{x})&=& \int_{\mathbb R^2}  |\mathbf{v}'| p(t,\mathbf{x},\mathbf{v}')\, d{\bf v}'.
\label{jescalar}
\end{eqnarray}
The variable coefficient $\alpha(c)$ is still defined by (\ref{Falpha}).
There are two reasons to replace $|\mathbf j|$ with $j$. From the
modeling point of view, the euclidean norm of $\mathbf j$ might
vanish under certain symmetry conditions, whereas in practice
the concentration of tumor angiogenic factor $c$ would still decrease
due to cell consumption. From the mathematical point of view,
$|\mathbf j(p)|= (j_1(p)^2+ j_2(p)^2)^{1/2} $ may bring about
lipschitzianity problems near zero.
This might cause uniqueness problems when $\mathbf j(p)$ 
approaches zero, which happens at infinity when $\mathbf x$
is allowed to vary in an unbounded domain. Our existence proof
holds for both choices, $|\mathbf j|$ and $j$. However, we can only
guarantee uniqueness in the latter case. }

{ Notice that we have also included a viscosity term 
$\Delta_{\mathbf x} p$. Adding a vanishing viscosity term
$\delta \Delta_{\mathbf x} p$, $\delta$ small, to problems with degenerate 
diffusion in that variable is a standard numerical strategy  to devise numerical
schemes ensuring positivity of solutions and avoiding sign related
artifacts. Our results extend to problems with asymmetric diffusion
$\sigma_1 \Delta_{\mathbf x} p + \sigma_2 \Delta_{\mathbf v} p$.
We have set $\sigma_1= \sigma_2 = \sigma$ for simplicity.}

The unknown $p$ represents a density of blood vessels. Positivity is therefore a key property of the solutions.  For positive $p$ the sign of the source term in (\ref{heat1}) cannot be controlled. A possibility to generate approximate solutions with a controlled sign is to freeze the integral coefficient and { to} include the linearized integral source in a linearized diffusion operator. The paper is organized as follows.
Section \ref{sec:linear} recalls and establishes results on existence of fundamental solutions,  solutions for initial value problems,  properties of such solutions and comparison principles, mostly  { for a linearized version of (\ref{heat1})-(\ref{heat2}) in which the dependence on $c$ is ignored}. Section \ref{sec:nonlinear0} constructs the unique nonnegative solution $p$
of the nonlinear problem when $\alpha=0$ as limit of solutions of an iterative scheme. Existence of nonnegative solutions for the scheme follows using fundamental solutions. Heat equations provide upper solutions yielding uniform $L^q$ bounds when $\alpha=0$. Energy inequalities produce the uniform bounds on derivatives required for compactness. The limiting function is the unique solution sought for. 
Once the strategy to handle the integral term is clear, section 
\ref{sec:coupling} considers the full coupled problem  (\ref{heat1})-(\ref{heat4}).  These results pave the way for the study of more realistic problems in which the heat operator is replaced by a Fokker-Planck 
operator with degenerate diffusion and transport terms, for which the 
theory of fundamental solutions is more involved 
\cite{carpio}. Our iterative schemes may be used for the numerical approximation of the solutions. The additional properties we establish on the solutions and the iterates may be exploited to render formal derivations of these models rigorous.


\section{Linear problem}
\label{sec:linear}

The key underlying linearized problems are:
\begin{eqnarray} Lp = \frac{\partial}{\partial t} p(t,\mathbf{x},\mathbf{v}) - \sigma
 \Delta_{\mathbf{x}\mathbf{v}} p(t,\mathbf{x},\mathbf{v}) + 
a(t,\mathbf{x}, \mathbf v) p(t,\mathbf{x},\mathbf{v}) 
= f(t,\mathbf{x},\mathbf{v}),  \label{heatlinear} \\
p(0,\mathbf{x},\mathbf{v})=p_0(\mathbf{x},\mathbf{v}), \label{heatdatum}
\end{eqnarray}
when $(\mathbf{x},\mathbf{v}) \in \mathbb{R}^2 \times \mathbb{R}^2$, $t \in [0,\infty)$, with
$a \in L^{\infty} ([0,\infty)\times \mathbb{R}^2\times \mathbb{R}^2)$, $\sigma \in \mathbb{R}^+$, 
$f \in L^{\infty}(0,\infty; L^{\infty}\cap L^{1} (\mathbb{R}^2 \times \mathbb{R}^2))$ and 
$p_0 \in L^{\infty}\cap L^{1} (\mathbb{R}^2 \times \mathbb{R}^2),$  as well as:
\begin{eqnarray} \frac{\partial}{\partial t} c(t,\mathbf{x}) - d
\Delta_{\mathbf{x}} c(t,\mathbf{x}) + a(t,\mathbf{x}) c(t,\mathbf{x}) 
= f(t,\mathbf{x}),  \label{heatlinearc} \\
c(0,\mathbf{x})=c_0(\mathbf{x}), \label{heatdatumc}
\end{eqnarray}
when $\mathbf{x} \in \mathbb{R}^2 $, $t \in [0,\infty)$, with
$a \in L^{\infty}  ([0,\infty)\times \mathbb{R}^2)$, 
$d \in \mathbb{R}^+$, 
$f \in L^{\infty}(0,\infty; L^q (\mathbb{R}^2))$ and 
$c_0 \in L^q (\mathbb{R}^2),$ $ 1 \leq q \leq \infty$.

Solutions may be constructed using fundamental solutions of the parabolic
operator, whose properties depend on the smoothness of the coefficient
$a$. We recall below the known theory of classical and weak fundamental
solutions and discuss additional bounds for solutions of (\ref{heatlinear})-(\ref{heatdatum})  and (\ref{heatlinearc})-(\ref{heatdatumc}).

\subsection{Existence using classical fundamental solutions}
\label{sec:linearclassical}

Existence of a classical smooth solution of (\ref{heatlinear})-(\ref{heatdatum}) follows from the theory of fundamental solutions for parabolic equations with smooth bounded coefficients \cite{friedman}. Let us assume that $a(t,\mathbf{x},\mathbf{v})$ 
is a continuous function satisfying
\begin{eqnarray}
|a(t,\mathbf{x}, \mathbf{v})-a(t,\mathbf{x}^0, \mathbf{v}^0)|\leq A 
|(\mathbf{x}-\mathbf{x}^0, \mathbf{v}-\mathbf{v}^0)|^{\beta}, 
\quad 0 < \beta < 1,
\label{holder}
\end{eqnarray}
for some $A>0$.  A fundamental solution of $Lu=0$ is a function
$\Gamma(t,\mathbf{x},\mathbf{v};\tau,\mathbf{x}',\mathbf{v}')$ defined for
$\mathbf{x},\mathbf{v},\mathbf{x}',\mathbf{v}' \in \mathbb{R}^2$ and $t,\tau \in \mathbb{R}^+$, $t>\tau$,
which verifies:
\begin{itemize}
\item[(i)] for fixed $(\tau,\mathbf{x}',\mathbf{v}')$, the equation $L \Gamma=0$ holds for all
$\mathbf{x},\mathbf{v}$, $t>\tau$,
\item[(ii)] for every continuous function $\psi$ 
\begin{eqnarray*}
\lim_{t\rightarrow \tau} \int_{\mathbb{R}^2} \int_{\mathbb{R}^2} \Gamma(t,\mathbf{x},\mathbf{v};\tau,\mathbf{x}',\mathbf{v}') \psi(\mathbf{x}',\mathbf{v}')
 d \mathbf{x}' d \mathbf{v}' = \psi(\mathbf{x},\mathbf{v}).
\end{eqnarray*}
\end{itemize}
According to { Theorems 10 and 11 in Chapter 1 of reference} \cite{friedman}, there exists a fundamental solution for our operator $L$ under  { hypothesis (\ref{holder})} on $a$, which satisfies the following bounds:
\begin{eqnarray}
0< \Gamma(t,\mathbf{x},\mathbf{v};\tau,\mathbf{x}',\mathbf{v}')| \leq C(T) (t-\tau)^{-n/2} e^{-\sigma^* {(|\mathbf{x}-\mathbf{x}'|^2+|\mathbf{v}-\mathbf{v}'|^2) \over 4(t-\tau)}}, \label{boundG}\\
|{\partial \Gamma(t,\mathbf{x},\mathbf{v};\tau,\mathbf{x}',\mathbf{v}')
\over \partial z_i}| \leq C(T) (t-\tau)^{-(n+1)/2} e^{-\sigma^* {(|\mathbf{x}-\mathbf{x}'|^2+|\mathbf{v}-\mathbf{v}'|^2) \over 4(t-\tau)}}, \label{bounddG}
\end{eqnarray}
for { $t \in [0,T]$,}
where $z_i=x_i$ or $v_i$, $i=1,2$, for $\sigma^* <\sigma$ and $C(T)>0$ { (see also reference \cite{porper}, pp. 124-125).} The { constants $C(T)$ appearing in estimates (\ref{boundG})-(\ref{bounddG})} depend on the parabolicity constant $\sigma$, the number of independent space and velocity variables $n$, the final time $T$, { the maximum modulus of the coefficients 
$M_0={\rm max}_{t\in[0,T], \mathbf x \in \mathbb{R}^2} |a(t, \mathbf x)|$ 
and the H\"older constants $A, \beta$ in inequality (\ref{holder})}, see reference \cite{porper}.
Given continuous functions $f(t,\mathbf{x},\mathbf{v})$ and $p_0(\mathbf{x},\mathbf{v})$ defined on $[0,\infty)\times \mathbb{R}^2 \times \mathbb{R}^2$ and $\mathbb{R}^2 \times \mathbb{R}^2$, respectively, and satisfying:
\begin{eqnarray*}
|f(t,\mathbf{x},\mathbf{v})| \leq  C_f e^{h (|\mathbf{x}|^2+|\mathbf{v}|^2)},
\quad
|p_0(\mathbf{x},\mathbf{v})| \leq  C_p e^{h (|\mathbf{x}|^2+|\mathbf{v}|^2)},
\end{eqnarray*}
for $C_f,C_p>0$ and $h <{\sigma \over 4 T}$, the function  
\begin{eqnarray}
p(t,\mathbf{x},\mathbf{v})= \int_{\mathbb{R}^2} \int_{\mathbb{R}^2} \Gamma(t,\mathbf{x},\mathbf{v};0,\mathbf{x}',\mathbf{v}') p_0(\mathbf{x}',\mathbf{v}') d \mathbf{x}' d \mathbf{v}' \nonumber \\
+ \int_0^t \int_{\mathbb{R}^2} \int_{\mathbb{R}^2} \Gamma(t,\mathbf{x},\mathbf{v};\tau,\mathbf{x}',\mathbf{v}') f(\tau,\mathbf{x}',\mathbf{v}')\, d\tau \, d \mathbf{x}' d \mathbf{v}'
\label{solint}
\end{eqnarray}
is a solution of problem (\ref{heatlinear})-(\ref{heatdatum}), as shown in reference \cite{friedman}. The positivity of the fundamental solution implies positivity of solutions for  $f \geq 0$ and $p_0 \geq 0$. The integral expression (\ref{solint}), together with the bound (\ref{boundG}), implies:
\begin{eqnarray}
\|p(t)\|_{\infty} \leq  C_{\infty} \left( \|p_0\|_{\infty} 
+ \int_0^t \|f(s)\|_{\infty} ds \right) \label{Linf} \\
\|p(t)\|_1 \leq  C_1 \left( \|p_0\|_1
+ \int_0^t \|f(s)\|_1 ds \right), \label{L1} 
\end{eqnarray}
where $C_{\infty}, C_1$ depend on  $\sigma$, $n$, $T$, $M_0$, $A$ and 
$\beta$. 

Bounds on classical fundamental solutions are studied more in detail in reference \cite{porper}, that admits $0 < \beta \leq 1$ and includes classical time derivatives and second order derivatives in space:
\begin{eqnarray}
|\partial_t^{m_0} \partial_{\mathbf{x}}^{m_x} \partial_{\mathbf{v}}^{m_v}
\Gamma(t,\mathbf{x},\mathbf{v};\tau,\mathbf{x}',\mathbf{v}')
| \leq {C(T) e^{-\sigma^* {(|\mathbf{x}-\mathbf{x}'|^2+|\mathbf{v}-\mathbf{v}'|^2) \over (t-\tau)}} \over (t-\tau)^{(n+2m_0+|m_x|+|m_v|)/2}} 
\label{boundsG}
\end{eqnarray}
for $0 < t-\tau \leq T$, $2m_0+|m_x|+|m_v| \leq 2$. The constants 
$\sigma^*$ and $C$ depend on $\sigma$, $n$, $M_0$, $A$, $\beta$ and 
$T$ {(see \cite{porper}, Theorem 1.1)}.

Existence results for measurable, bounded or integrable data follow by a regularization procedure. We smooth the data using convenient mollifiers \cite{brezis,lions}: $p_0= \rho_{\varepsilon} * p_0$ and  $f= \psi_{\varepsilon} * f$ in adequate variables. These $C^{\infty}$ families converge to $p_0$ and $f$ in $L^q$, $1 \leq q < \infty$, as $\varepsilon$ tends to zero, and are bounded when $p_0, f \in L^q$. When $q=\infty$, we have weak* convergence. The corresponding classical sequences of solutions, and their derivatives, are bounded in  $L^q$. Passing to the limit in the linear equation, the limit is a weak solution. Passing to the limit in the integral equation, it verifies a similar integral equation. The limit solution inherits nonnegativity for positive data. It also inherits the $L^q$ bounds (\ref{Linf})-(\ref{L1}).

\subsection{Existence using weak fundamental solutions}
\label{sec:linearweak}

When the coefficient $a$ is only measurable and bounded, existence of weak fundamental solutions has been established. A measurable function of the form $\Gamma(t,\mathbf{x},\mathbf{v};\tau,\mathbf{x}',\mathbf{v}')$ is a weak fundamental solution of the initial value problem for (\ref{heatlinear})-(\ref{heatdatum}) if the function 
\begin{eqnarray}
P_{t,\tau}\psi (\mathbf{x},\mathbf{v})=\int_{\mathbb{R}^2 \times \mathbb{R}^2} 
 \Gamma(t,\mathbf{x},\mathbf{v}; \tau,\mathbf{x}',\mathbf{v}')
\psi({\mathbf x}',{\mathbf v}') d{\mathbf x}' d{\mathbf v}'
\label{integralweak}
\end{eqnarray}
satisfies:
\begin{eqnarray}
{\partial \over \partial t} P_{t,\tau}\psi (\mathbf{x},\mathbf{v}) =
[\sigma \Delta_{\mathbf{x},\mathbf{v}} - a(t,\mathbf{x},\mathbf{v}) ] P_{t,\tau}\psi (\mathbf{x},\mathbf{v}) \label{weak1} \\
{\rm lim}_{t\rightarrow \tau} P_{t,\tau}\psi(\mathbf{x},\mathbf{v}) =
\psi({\mathbf x},{\mathbf v}) \label{weak2}
\end{eqnarray}
for any continuous function $\psi$ with compact support, $\mathbf{x},\mathbf{v},
\mathbf{x}',\mathbf{v}' \in \mathbb{R}^2 \times \mathbb{R}^2$, $t,\tau \in [0,\infty)$, $t>\tau$.
The weak solution of  the initial value problem (\ref{weak1})-(\ref{weak2}) is unique since $\sigma$ is constant \cite{kusuoka}.
A weak solution $P_{t,\tau}{\psi}(\mathbf{x},\mathbf{v})$ of
(\ref{weak1})-(\ref{weak2})  satisfies
\begin{eqnarray*}
 \int_0^t \hskip -2mm \int_{\mathbb{R}^2 \times \mathbb{R}^2} 
\hskip -9mm P_{s,0}\psi (\mathbf{x},\!\mathbf{v}) 
{
[{\partial \over \partial t} \!+\! \sigma \Delta_{\mathbf{x}\mathbf{v}}  
\!\!-\! a(s,\!\mathbf{x},\!\mathbf{v}\!) ] 
\phi(t,\!{\mathbf x},\!{\mathbf v})} d{\mathbf x} d{\mathbf v} ds \!+\!\!\!
 \int_{\mathbb{R}^2 \times \mathbb{R}^2} \hskip -8mm \psi({\mathbf x},\!{\mathbf v})
\phi(0,\!{\mathbf x},\!{\mathbf v}) d{\mathbf x} d{\mathbf v} \!=\!0,
\end{eqnarray*}
{ for
any $\phi \in C^{\infty}_c([0,T) \times \mathbb{R}^2 \times \mathbb{R}^2)$}.  

Under hypothesis on the second order operator that are satisfied for constant positive 
$\sigma$ and measurable, bounded $a$,  a fundamental solution 
$\Gamma$ satisfying the bounds \cite{aronson1,aronson2,norris,kusuoka}:
\begin{eqnarray}
{C_1 e^{-C_1 (t-\tau)} \over (t-\tau)^{n/2}}
e^{-\gamma_1 {(|\mathbf{x}-\mathbf{x}'|^2+|\mathbf{v}-\mathbf{v}'|^2) \over t-\tau}} \leq 
\Gamma(t,\mathbf{x},\mathbf{v};\tau,\mathbf{x}',\mathbf{v}')   \label{lowerboundG} \\
\Gamma(t,\mathbf{x},\mathbf{v};\tau,\mathbf{x}',\mathbf{v}') 
\leq {C_2 e^{C_2 (t-\tau)} \over
(t-\tau)^{n/2}} e^{-\gamma_2 {(|\mathbf{x}-\mathbf{x}'|^2+|\mathbf{v}-\mathbf{v}'|^2) 
\over t-\tau}},
\label{upperboundG}
\end{eqnarray}
exists for $t, \tau \in [0, \infty)$ such that $\tau<t$ and $\mathbf{x},\mathbf{v},\mathbf{x}',\mathbf{v}' \in \mathbb{R}^2$. The dimension $n=4$ in our particular case. The constants $C_1, C_2, \gamma_1, \gamma_2 $ depend on $\sigma$ and  $\|a\|_{L^{\infty}}$. 

The existence and regularity of weak fundamental solutions for {  general para-\ bolic} problems of the form
$u_t=\nabla \cdot (\sigma \nabla u) + b \cdot \nabla u - a u$ with measurable coefficients 
has been studied in a series of papers. $\beta$ - H\"older continuity, $\beta \in (0,1],$ is discussed in \cite{aronson2,stroock,kusuoka}.  
When $b$ is not continuous, differentiability cannot be expected. 
Depending on the regularity of the coefficients \cite{norris,kusuoka}, 
{ the equivalent of initial value problem (\ref{weak1})-(\ref{weak2}) 
for general parabolic  operators}  should be understood in a merely weak sense.

In our case, $\sigma$ is constant, $b=0$, and $a$ is bounded. Therefore, 
both the fundamental solution $\Gamma$ and $P_{t,\tau} \psi$ are 
smoother. Once its existence is guaranteed, the fundamental solution 
can be seen as a solution of a heat equation with a source 
$- a p=-a \Gamma$, bounded in terms of a heat kernel. {
It admits the integral expression
\begin{eqnarray}
\Gamma(t, {\mathbf x}, {\mathbf v};  \tau, {\mathbf x}', {\mathbf v}')  = 
G(t - \tau, {\mathbf x} - {\mathbf x}', \mathbf{v} - \mathbf{v}')  \nonumber \\
- \int_{\tau}^t \int_{\mathbb R^2 \times \mathbb R^2}  
G(t-s, {\mathbf x}-{\boldsymbol \xi}, \mathbf{v}-{\boldsymbol \nu})  
a(s, {\boldsymbol \xi}, {\boldsymbol \nu})
\Gamma(s,  {\boldsymbol \xi}, {\boldsymbol \nu};   
\tau,  {\mathbf x}', \mathbf{v}') 
d{\boldsymbol \xi} d{\boldsymbol \nu} ds,
\label{eq:integralGa}
\end{eqnarray}
where $G$ is the heat kernel for diffusivity $\sigma$. 
In fact, $\Gamma$ can be constructed as the solution of integral equation 
(\ref{eq:integralGa}) using an iterative scheme that yields bound 
(\ref{boundG}), but with coefficients depending on $\|a\|_{\infty}$, $n$, 
$\sigma$ and $T$, for $t \in[0,T]$. This way of reasoning is standard in 
kinetic models, see  \cite{victoryclassical,carpio} and references therein.
The derivatives of $\Gamma$ satisfy a similar integral expression. 
For $z_i=x_i$ or $z_i=v_i$, $i=1,2$:
\begin{eqnarray}
{\partial \over \partial z_i}\Gamma(t, {\mathbf x}, {\mathbf v};  
\tau, {\mathbf x}', {\mathbf v}')  = {\partial \over \partial z_i} 
G(t - \tau, {\mathbf x} - {\mathbf x}', \mathbf{v} - \mathbf{v}')  \nonumber \\
- \int_{\tau}^t \int_{\mathbb R^2 \times \mathbb R^2}  
{\partial \over \partial z_i} 
G(t-s, {\mathbf x}-{\boldsymbol \xi}, \mathbf{v}-{\boldsymbol \nu})  
a(s, {\boldsymbol \xi}, {\boldsymbol \nu})
\Gamma(s,  {\boldsymbol \xi}, {\boldsymbol \nu};   
\tau,  {\mathbf x}', \mathbf{v}') 
d{\boldsymbol \xi} d{\boldsymbol \nu} ds.
\label{eq:integraldGa}
\end{eqnarray}
The coefficient $a$ being bounded, inequality (\ref{upperboundG}) 
allows us to obtain estimates of the form (\ref{bounddG}), like 
in the classical case, but with coefficients depending on $\|a\|_{\infty}$, $n$, 
$\sigma$ and $T$, for $t \in[0,T]$.}  

{
We summarize these observations for the problems we consider here
in the following Lemmas. \\

{\bf Lemma 2.1.} {\it When $a$ is a bounded function, the initial value 
problem (\ref{heatlinear})-(\ref{heatdatum}) has a fundamental solution 
$\Gamma$ satisfying estimates (\ref{lowerboundG})-(\ref{upperboundG}) 
and (\ref{bounddG}) with parameters depending on the norm $\|a\|_{\infty}$, 
the dimension $n$, the diffusivity $\sigma$ and $T>0$, for $t \in[0,T]$.
An analogous result holds for (\ref{heatlinearc})-(\ref{heatdatumc}).} \\

{\bf Proposition 2.2.} {\it  For any 
$a \in  L^\infty ([0,\infty) \! \times \! \mathbb{R}^2 \! \times \! \mathbb{R}^2)$, 
$p_0 \in L^\infty \cap L^1(\mathbb{R}^2  \! \times  \! \mathbb{R}^2)$ and  
$f \!  \in \!  L^\infty (0,T;L^\infty  \cap L^1 (\mathbb{R}^2 \!  \times \! \mathbb{R}^2))$,  there exists a unique solution 
$p \! \in \!  C([0,\! T];L^\infty \cap L^1 (\mathbb{R}^2 \! \times \! \mathbb{R}^2))$  of the initial value problem (\ref{heatlinear})-(\ref{heatdatum})
satisfying the weak formulation:
\begin{eqnarray}
 \!\! \int_0^t \hskip -2mm \int_{\mathbb{R}^2 \times \mathbb{R}^2} 
 \hskip -3mm
p  [ {\partial \phi \over \partial t}  \! + \!
\sigma \Delta_{\mathbf{x}\mathbf{v}} \phi  \! - \! a\phi] 
d{\mathbf x} d{\mathbf v} ds
+ \!\! \int_0^t \hskip -2mm \int_{\mathbb{R}^2 \times \mathbb{R}^2} 
 \hskip -8mm f \phi  d{\mathbf x} d{\mathbf v} ds 
+  \hskip -2mm \int_{\mathbb{R}^2 \times \mathbb{R}^2} 
\hskip -8mm p_0 \phi(0) d{\mathbf x} d{\mathbf v} = 0, 
\label{fundweak}
\end{eqnarray}
for $\phi  \in C^{\infty}_c([0,T) \times \mathbb{R}^2 \times \mathbb{R}^2).$
This solution admits the integral
expression (\ref{solint}) in terms of the fundamental solution $\Gamma$
given by (\ref{eq:integralGa}). It is positive for 
positive $f$ and $p_0$ and satisfies estimates (\ref{Linf})-(\ref{L1}) 
with constants depending on $\sigma$, $T$, $\|a\|_{\infty}$ and
the dimension. The solution has the same regularity as solutions
of heat equations with $L^q$ data and sources, $1\leq q \leq \infty$.}  

{\bf Proof.} We can either exploit (\ref{integralweak}) and the semigroup 
theory \cite{pazy}, or construct the solution as limit of classical solutions 
as in reference \cite{carpio}. 
Following \cite{carpio}, let us consider regularized coefficients 
$a_{k} = a * \phi_k$ where $\phi_k$ is a positive $C^{\infty}$
mollifying family (see reference \cite{brezis}, p. 108). We  then have
$\| a_kÊ\|_\infty \leq \|aÊ\|_\infty$ and $a_k \rightharpoonup a$ 
in $L^\infty_{t\mathbf x \mathbf v}$ weak* when
$k \rightarrow \infty$ (see reference \cite{brezis},  p. 126). 
Let us first assume that $p_0,f \in C^\infty_c$ are smooth functions with 
compact support. Let $p_k$ be the corresponding classical solutions of 
(\ref{heatlinear})-(\ref{heatdatum}) with coefficient $a_k$,
satisfying the integral equations (\ref{solint}) in terms
of fundamental solutions $\Gamma_k$.
Thanks to estimates (\ref{lowerboundG})-(\ref{upperboundG})
and the fact that $\| a_kÊ\|_\infty$ is uniformly bounded,
$p_k$ are bounded in $L^2(0,T;L^2_{\mathbf x \mathbf v})$.
Then, the energy inequality (see \cite{cazenave}, also Lemma 2.9 below)
provides a uniform bound on
${\partial p_k \over \partial z_i}$  
in $L^2(0,T;L^2_{\mathbf x \mathbf v})$ for 
$z_i=v_i$ or $x_i$, $i=1,2$.   
Therefore, $p_k$ is bounded in $L^2(0,T;H^1_{\mathbf x \mathbf v})$. 
Using equation (\ref{heatlinear}), the time derivative ${\partial
p_k \over \partial t}$ is bounded in $L^\infty(0,T;H^{-1}_{\mathbf x \mathbf v})$ (see \cite{brezis}, proposition 9.20).
Since the injection of $H^{1}(\Omega)$ in $L^2(\Omega)$ is
compact for any smooth bounded $\Omega$ (see \cite{brezis}, theorem 9.16), classical compactness results 
(see Theorem 12.1 in \cite{lions}, Corollary 4 in \cite{simon}) 
imply that $p_k$ is 
compact in $L^2(0,T;L^2_{\mathbf x \mathbf v}(\Omega))$, for any bounded
$\Omega$.  We can extract subsequences converging strongly and
pointwise to $p$ in any $\Omega$. 
Choosing $\Omega=B(0,M)$, for integers $M$ tending to infinity, a diagonal 
extraction procedure allows us to
obtain a subsequence $p_{k'}$ converging pointwise to a limit $p$
in $[0,T] \times \mathbb{R}^2 \times \mathbb{R}^2$, strongly in 
$L^2(0,T;L^2_{loc})$, and  weakly in $L^2(0,T;L^2_{\mathbf x \mathbf v})$.
For $\phi  \in C^{\infty}_c([0,T) \times \mathbb{R}^2 \times \mathbb{R}^2),$
the weak formulation of problem (\ref{heatlinear})-(\ref{heatdatum})
with coefficient $a_k$ reads:
\begin{eqnarray}
 \!\! \int_0^t \hskip -2mm \int_{\mathbb{R}^2 \times \mathbb{R}^2} 
 \hskip -5mm
p_{k}  [ {\partial \phi \over \partial t}  \! + \!
\sigma \Delta_{\mathbf{x}\mathbf{v}} \phi  \! - \! a_{k} \phi] 
d{\mathbf x} d{\mathbf v} ds
+ \!\! \int_0^t \hskip -2mm \int_{\mathbb{R}^2 \times \mathbb{R}^2} 
 \hskip -8mm f \phi  d{\mathbf x} d{\mathbf v} ds 
+  \hskip -2mm \int_{\mathbb{R}^2 \times \mathbb{R}^2} 
\hskip -8mm p_0 \phi(0) d{\mathbf x} d{\mathbf v} = 0.
\label{fundweakk}
\end{eqnarray}
The support of $\phi$ being a compact set contained in 
$[0,T]\times \Omega$,
$\Omega$ bounded, we may select a subsequence
converging to $p$ in $L^2(0,T;L^2(\Omega))$,  and therefore in 
$L^1(0,T;L^1(\Omega))$. Taking limits, $p$ is a weak solution of (\ref{heatlinear})-(\ref{heatdatum}) with coefficient $a$.

Now, let us pass to the limit in the integral expressions (\ref{solint})
for $p_{k'}$ and $\Gamma_{k'}$. Thanks to inequality (\ref{upperboundG}), 
the fundamental solutions $\Gamma_{k'}$ are bounded in 
$L^r_{t \mathbf x \mathbf v}((0,T)\times I\!\!R^2 \times I\!\!R^2)$ 
for any $ r\in (1, n/(n-2))$, $n=4$. 
%
%
Therefore, a subsequence converges weakly in $L^r_{t \mathbf x
\mathbf v}$ to a limit $\Gamma$. Taking limits in  the integral 
expressions, we see that the solution $p$ of 
(\ref{heatlinear})-(\ref{heatdatum}) with coefficient $a$ verifies (\ref{solint}).
Setting $f=0$, $\Gamma$ fulfills the definition of fundamental solution 
for $a$. 

Once identity (\ref{solint}) is established for solutions of
(\ref{heatlinear})-(\ref{heatdatum})  with bounded coefficient $a$
and $C^\infty_c$ data, we extend it to $L^q$ data by density.
When $1 \leq q < \infty$, there are sequences
$p_{k,0} \in C^\infty_c$ and $f_k \in C^\infty_c$ 
converging to $p_0$ in $L^q_{\mathbf x \mathbf v}$ and 
$f$ in $L^1_t L^q_{\mathbf x \mathbf v}$ (see \cite{brezis}, corollary 4.23).
Identity (\ref{solint}) and inequality (\ref{upperboundG}) 
yield a uniform bound on $\|p_{k}(t)\|_r$ for $t \in [0,T]$, $r \geq q$. 
Therefore, we can extract a subsequence converging weakly to a 
limit $p$. Taking limits in the weak formulation and the integral
equation for $p_k$, we prove (\ref{solint}) for a weak
solution $p$ of the initial value problem (\ref{heatlinear})-(\ref{heatdatum}) 
with data $p_0$ and $f$. 

%
%

Using identity (\ref{solint}), the positivity of the solution follows from 
the positivity of the data.
The upper bound (\ref{upperboundG}) on the fundamental solution
yields estimates (\ref{Linf})-(\ref{L1}) with constants depending on 
$\sigma$, $n$, $T$, $\|a\|_{\infty}$. 

To prove uniqueness, assume we have two solutions $p_1,p_2 \in 
L^\infty(0,T;L^q_{\mathbf x \mathbf v})$ of (\ref{heatlinear})-(\ref{heatdatum}) 
with initial and source data in $L^q$ and coefficient $a$ bounded. 
The difference $\overline{p}=p_1-p_2$ is a solution of a heat equation
with $\overline{p}(0)=0$ and source $-a \overline{p}$. The integral
expression for solutions of heat equations yields:
\begin{eqnarray*}
\overline{p}(t)= - \int_0^t  G(t-s)*a(s)\overline{p}  ds
\Rightarrow \|\overline{p}(t) \|_q \leq \|a\|_{\infty} \int_0^t
\|\overline{p}(s) \|_q  ds.
\end{eqnarray*}
This Gronwall inequality implies $\|\overline{p} (t)\|_q =0$ for any
$t \in[0,T]$. \\

Similar results hold suppressing the dependence on the variable 
$\mathbf v$.  We state the result for data in $L^q_{\mathbf x}$ spaces 
since we do not assume 
$c_0 \in L^\infty_{\mathbf x} \cap L^1_{\mathbf x}$.

{\bf Proposition 2.3.} {\it  For any 
$a \in   L^\infty ([0,\infty) \times \mathbb{R}^2)$, 
$c_0 \in L^q(\mathbb{R}^2)$,  and 
$f \in L^\infty(0,T;L^q (\mathbb{R}^2))$, $1 \leq q \leq \infty$,
there exists a unique solution $c \in C([0,T];L^\infty (\mathbb{R}^2))$ 
of the initial value problem 
(\ref{heatlinearc})-(\ref{heatdatumc}). This solution admits the integral
expression (\ref{solint}) in terms of the fundamental solution $\Gamma$
(suppressing the dependence on $\mathbf v$).
The solution of (\ref{heatlinearc})-(\ref{heatdatumc}) is positive for 
positive $f$ and $p_0$ and satisfies estimates (\ref{Linf})-(\ref{L1}) 
with constants depending on $\sigma$, $T$, $\|a\|_{\infty}$ and
the dimension.}

{\bf Proof.} The proof is similar to that of Proposition 2.2, except
for the uniform bounds on $p_k$ and ${\partial p_k \over \partial z_i}$. 

Let us first assume that $q=\infty$.
Thanks to estimates (\ref{lowerboundG})-(\ref{upperboundG})
and the fact that $\| a_kÊ\|_\infty$ is uniformly bounded,
$p_k$ are bounded in $L^\infty(0,T;L^\infty_{\mathbf x \mathbf v})$.
Moreover, differentiating (\ref{solint}) and using the estimates
on the derivatives of $\Gamma_k$ referred to in Lemma 2.1,
we see that ${\partial p_k \over \partial z_i}$ are uniformly
bounded in $L^{r_1}(0,T;L^\infty_{\mathbf x \mathbf v})$ for 
$z_i=v_i$ or $x_i$, $i=1,2$, and $r_1\in [1,2)$. Therefore,
$p_k$ is uniformly bounded in 
$L^{r_1}(0,T;W^{1,r_2}_{\mathbf x \mathbf v}(\Omega))$ 
for any $r_2 \in [1,\infty)$ and any $\Omega$ bounded.
Using equation (\ref{heatlinear}), the time derivative ${\partial
p_k \over \partial t}$ is bounded in 
$L^{r_1}(0,T;W^{-1,r_2'}_{\mathbf x \mathbf v})$ when $r_2<\infty$ 
(see \cite{brezis}, proposition 9.20).
Since the injection of $W^{1,r_2}(\Omega)$ in $L^{r_2}(\Omega)$ is
compact for any smooth bounded 
$\Omega$ (\cite{brezis}, theorem 9.16), classical compactness results \cite{aubin,lions,simon}  imply that $p_k$ is compact 
in $L^{r_1}(0,T;L^{r_2}_{\mathbf x \mathbf v}(\Omega))$.
As in Proposition 2.2, a diagonal  extraction procedure allows us to
obtain a subsequence $p_{k'}$ converging pointwise to a limit $p$
in $[0,T] \times \mathbb{R}^2 \times \mathbb{R}^2$, strongly in 
$L^{r_1}(0,T;L^{r_2}_{loc})$, and  weak* in 
$L^\infty(0,T;L^\infty_{\mathbf x \mathbf v})$. In particular,
we have strong convergence in $L^1(0,T;L^1_{loc})$ which
allows us to pass to the limit in the weak formulation of the
equations (\ref{fundweakk}) and establish that $p$ satisfies  (\ref{fundweak}).
When $q<\infty$, a similar proof works.

Identity (\ref{solint}) when $q< \infty$ is proven as in Proposition
2.2. When $q=\infty$, we first have to extend it
to solutions with coefficient $a$ and differentiable bounded data,
with bounded derivatives. The proof proceeds as the proof in 
Proposition 2.2 for $C^\infty_c$
data since classical solutions for these data satisfy the integral
equation, and bounds in $L^\infty(0,T;W^{1,\infty}_{\mathbf x \mathbf v})$
imply bounds in $L^2(0,T;H^{1}(\Omega)),$ $\Omega$ bounded.
This allows to pass to the limit in the weak formulations of the
initial values problem with coefficient $a_k$. Also, the
associated fundamental solutions $\Gamma_k$ satisfy the Dunford-Pettis
criterion for weak compactness in $L^1$ (see reference \cite{brezis},
theorem 4.30), which allows us taking limits in the integral equation.
Once  (\ref{solint}) is established for solutions with coefficient $a$ and 
differentiable bounded data, with bounded derivatives, $L^\infty$ data
are handled approximating $p_0$ and $f$ by mollified
sequences $p_{k,0}$ and $f_k$ tending to $p_0$ and $f$ in $L^\infty$
weak*. 
}

\subsection{Comparison principle and heat estimates}
\label{sec:comparison}

This section recalls comparison principles and basic $L^r-L^q$ estimates
for solutions of diffusion problems. \\

{
{\bf Lemma 2.4.} { \it Let $p^{(1)}$, $p^{(2)}$ be the solutions of
the initial value problem
(\ref{heatlinear})-(\ref{heatdatum}) with bounded coefficient
$a $ and data 
$f^{(1)}$, $p_{0}^{(1)}$ and  $f^{(2)}$, $p_{0}^{(2)}$,
respectively, constructed in Proposition 2.2. Assume that:
\begin{eqnarray}
f^{(1)} \leq f^{(2)}, \quad p_{0}^{(1)} \leq p_{0}^{(2)}. \label{comp1}
\end{eqnarray}
Then, the corresponding solutions $p^{(1)}$ and $p^{(2)}$ preserve 
the ordering:
\begin{eqnarray}
p^{(1)} \leq p^{(2)}.\label{comp2}
\end{eqnarray}
In particular, any solution $p$ is nonnegative if $p_0\geq 0$ and
$f \geq 0$.
The same positivity and comparison principles hold for 
(\ref{heatlinearc})-(\ref{heatdatumc}).}
} \\
{\bf Proof.}
It follows from the positivity of the fundamental solutions and the integral 
expression (\ref{solint}) for the solutions of the initial value problem 
(\ref{heatlinear})-(\ref{heatdatum}). Similarly for 
(\ref{heatlinearc})-(\ref{heatdatumc}).\\

{ \bf Lemma 2.5. {\it When the bounded coefficient
$a \geq 0$, any positive solution 
$p$ of the initial value problem (\ref{heatlinear})-(\ref{heatdatum}) 
with $L^q$ data is 
bounded from above by a solution of a heat equation with the same 
initial and source data. Moreover, the following estimates hold for
any $q \in[1, \infty]$: 
\begin{eqnarray} 
\| p \|_q &\leq& \| p_0 \|_q + t \,{\rm max}_{s\in [0,t]} \|f(s)\|_q, 
\label{LqLrp0f0} \\
\| p \|_r &\leq& C_1 t^{-({1\over q}-{1 \over r}){n \over 2}} \| p_0 \|_q
+ C_2 t^{-({1\over q}-{1 \over r}){n \over 2} +1} 
{\rm max}_{s\in [0,t]} \|f(s)\|_q, 
\label{LqLrp0f}
\end{eqnarray}
provided $r\geq q$, $({1\over q}-{1 \over r}){n \over 2} <1$, $n$ being
the dimension. Analogous estimates hold for solutions of
(\ref{heatlinearc})-(\ref{heatdatumc}) adapting the dimension.}Ê}
\\
{\bf Proof.}
Notice that $p$ is the solution of the heat equation with source 
$g=f - a p \leq f$. Let $u$ be the solution of:
\begin{eqnarray} \frac{\partial}{\partial t} u(t,\mathbf{x},\mathbf{v}) - \sigma
 \Delta_{\mathbf{x}\mathbf{v}} u(t,\mathbf{x},\mathbf{v}) 
= f(t,\mathbf{x},\mathbf{v}), \quad
u(0,\mathbf{x},\mathbf{v})=p_0(\mathbf{x},\mathbf{v}). \label{heat}
\end{eqnarray}
The solution of problem (\ref{heat}) admits
integral expressions in terms of the heat kernel $G(t,\mathbf{x},\mathbf{v})$.
It is then straightforward that:
\begin{eqnarray}
p(t)= G(t)*p_0 + \int_0^t  G(t-\tau)*[f(\tau)-a(\tau)p(\tau)]  d\tau  
\nonumber \\
\leq u(t)= G(t)*p_0 + \int_0^t  G(t-\tau)*f(\tau)  d\tau, 
\label{heatint}
\end{eqnarray}
where $*$ denotes convolution in the $\mathbf{x},\mathbf{v}$ variables.
Setting $f=0$, the well known $L^r-L^q$ estimates for heat operators
$\| u \|_q = \| G(t)*p_0 \|_q$ follow \cite{cazenave}:
\begin{eqnarray} 
\| u \|_q &\leq& \|G(t) \|_1 \| p_0 \|_q \leq \| p_0 \|_q, \\
\| u \|_r &\leq& \|G(t) \|_{q'} \| p_0 \|_q \leq C_{q'} t^{-({1\over q}-{1 \over r}){n \over 2}}
\| p_0 \|_q,  \; 1/r= 1/q+1/q'-1,
\label{LqLrp0}
\end{eqnarray}
for $r \geq q.$ The parameters $n$ stands for the space-velocity dimension.
When $f \neq 0$ we find estimates (\ref{LqLrp0f0})-(\ref{LqLrp0f}) for $u$.
They extend to $p$ since $p \leq u$.  Analogous arguments work
for (\ref{heatlinearc})-(\ref{heatdatumc}), with $n$ representing
only the spatial dimension.


\subsection{Velocity decay}
\label{sec:velocity}

{
Integral expressions for solutions of the initial value problem
(\ref{heatlinear})-(\ref{heatdatum}) yield additional information
on their velocity decay, depending on the initial and
source data. \\

{\bf Lemma 2.6.} {\it For the solution constructed in Proposition 2.2,
let us assume that
$|\mathbf v|^{\beta} f \in L^\infty(0,T;L^1(\mathbb R^2 \times \mathbb R^2))$
and
$|\mathbf v|^{\beta} p_0 \in L^1( \mathbb R^2 \times \mathbb R^2)$,
for $\beta >0$. Then, the solution $p$ of
(\ref{heatlinear})-(\ref{heatdatum}) satisfies
$|\mathbf v|^{\beta} p \in L^\infty(0,T;L^1(\mathbb R^2 \times \mathbb R^2)).$

Moreover, if 
$|\mathbf v|^\beta p_0 \in L^q_{\mathbf x}(\mathbb R^2;L^1_{\mathbf v}(\mathbb R^2))$
and
$ |\mathbf v|^\beta f \in L^\infty(0,T;L^q_{\mathbf x}(\mathbb R^2;L^1_{\mathbf v}(\mathbb R^2))),$ for
$ 1  \leq q \leq \infty$, $\beta=0,1,2$, then, the solution $p$ of
the problem (\ref{heatlinear})-(\ref{heatdatum}) satisfies
$|\mathbf v|^\beta p \in L^\infty(0,T;L^q_{\mathbf x}(\mathbb R^2;L^1_{\mathbf v}(\mathbb R^2)))$ for $\beta=0,1,2$.
}

{\bf Proof.}
Using  the integral expression (\ref{solint}) for $p$,
taking absolute values, multiplying by $|\mathbf v|^\beta$ 
and integrating we obtain:
\begin{eqnarray} 
\int_{I\!\!R^2 \times I\!\!R^2}  \hskip -10mm
|\mathbf v|^\beta |p(t, \! {\mathbf x}, \! {\mathbf v})| 
d \mathbf v d \mathbf x \leq  \hskip -2mm
\int_{I\!\!R^2\times I\!\!R^2\times I\!\!R^2 \times I\!\!R^2}  
\hskip - 2.2cm |\mathbf v|^\beta
\Gamma (t,\! {\mathbf x},\! {\mathbf v};  0, \!{\mathbf x}', \!{\mathbf v}')  
|p_0( \!{\mathbf x}',\! {\mathbf v}')| d{\mathbf x}' d{\mathbf v}' 
 d \mathbf x d \mathbf v + \nonumber \\
\int_0^t \hskip -2mm
\int_{I\!\!R^2\times I\!\!R^2\times I\!\!R^2 \times I\!\!R^2}  
\hskip - 2.2cm |\mathbf v|^\beta
\Gamma (t,\! {\mathbf x},\! {\mathbf v};  \!s, \!{\mathbf x}', \!{\mathbf v}')  
|f(s, \!{\mathbf x}',\! {\mathbf v}')| d{\mathbf x}' d{\mathbf v}' 
d s d \mathbf x d \mathbf v,
\end{eqnarray}
for any $t \in[0,T]$.
Thanks to estimate (\ref{upperboundG}),  the first integral is 
bounded from above in terms of a heat kernel $G$:
\begin{eqnarray}
C_1(T)  t^{\beta\over 2} \hskip -2mm 
\int_{\mathbb R^2 \times \mathbb R^2 \times \mathbb R^2 \times \mathbb R^2}  
\hskip -4mm
{|\mathbf v - \mathbf v'|^\beta \over t^{\beta/2}}
G(t,{\mathbf x}-{\mathbf x}',{\mathbf v}-{\mathbf v}') 
|p_0( {\mathbf x}', {\mathbf v}')| d{\mathbf x}' d{\mathbf v}' 
d \mathbf x d \mathbf v 
\nonumber \\
+ C_1(T) 
\int_{\mathbb R^2 \times \mathbb R^2 \times \mathbb R^2 \times \mathbb R^2}   
\hskip -4mm
G(t,{\mathbf x}-{\mathbf x}',{\mathbf v}-{\mathbf v}')  
|\mathbf v'|^\beta |p_0({\mathbf x}', {\mathbf v}')| d{\mathbf x}' 
d{\mathbf v}' d \mathbf x d \mathbf v Ê\nonumber  \\
\leq 
C_2(T) \|Êp_0\|_{L^1_{\mathbf x \mathbf v}}
+ C_3(T) \|Ê|Ê\mathbf v|^\beta p_0\|_{L^1_{\mathbf x \mathbf v}}.
\label{momentoj}
\end{eqnarray}
Proceeding in a similar way,
the integral involving $f$ is bounded  by
$\hat C_2(T) \|Êf\|_{L^\infty_t L^1_{\mathbf x \mathbf v}}$ $
+ \hat C_3(T) \|Ê|Ê\mathbf v|^\beta f\|_{L^\infty_t L^1_{\mathbf x \mathbf v}}.$

For the last part, integrating (\ref{solint})
with respect to $\mathbf v$, using (\ref{upperboundG}),
and 
\[ \int_{\mathbb R^2} d \mathbf v 
(t-\tau)^{-4/2} e^{-\gamma_2 {(|\mathbf{x}-\mathbf{x}'|^2+
|\mathbf{v}-\mathbf{v}'|^2) \over t-\tau} } = C_1 (t-\tau)^{-2/2}
e^{-\gamma_2 {|\mathbf{x}-\mathbf{x}'|^2 \over t-\tau} }
= K(t-s,{\mathbf x}-{\mathbf x}'), \]
we see that:
\begin{eqnarray}
\int_{\mathbb R^2} \hskip -3mm p d \mathbf v \leq  C_2 
\hskip -2mm \int_{\mathbb R^2\times \mathbb R^2} 
\hskip -9mm K(t,\!{\mathbf x}\!-\!{\mathbf x}'\!) 
p_0( {\mathbf x}'\!,\! {\mathbf v}'\!) d{\mathbf x}' d{\mathbf v}' \!+\!
C_2 \hskip -2mm \int_0^t \hskip -2mm
\int_{\mathbb R^2 \times \mathbb R^2} 
\hskip -9mm K(t\!-\!s,\!{\mathbf x}\!-\!{\mathbf x}'\!) 
f(s, \!{\mathbf x}' \!, \!{\mathbf v}'\!) d{\mathbf x}' \!d{\mathbf v}' \!ds.
\end{eqnarray}
This yields the $L^q_{\mathbf x}$ estimate on 
$\int_{\mathbb R^2}  p d \mathbf v$. 

For $\int_{\mathbb R^2} |\mathbf v| p d \mathbf v$, we multiply
identity (\ref{solint}) by $|\mathbf v|$, integrate with respect
to $\mathbf v$,  replace in the right hand side
$|\mathbf v|$ with $ |\mathbf v - \mathbf v'| + |\mathbf v'|$
and use estimate (\ref{upperboundG}).
We then notice that
\[ \int_{\mathbb R^2}  d \mathbf v 
{|\mathbf v -\mathbf v'|Ê\over (t-\tau)^{2}} 
e^{-\gamma_2 {(|\mathbf{x}-\mathbf{x}'|^2+
|\mathbf{v}-\mathbf{v}'|^2) \over t-\tau} } = {\hat C_1 \over t-\tau}
e^{-\gamma_2 {|\mathbf{x}-\mathbf{x}'|^2 \over t-\tau} }
= \tilde K(t-s,{\mathbf x}-{\mathbf x}'). \]
We finally find:
\begin{eqnarray}
\int_{\mathbb R^2}  |\mathbf v| p d \mathbf v \leq  \hat C_2 
\Big[ K(t)*  \int_{\mathbb R^2} |\mathbf v' | p_0 d \mathbf v' +
 \int_0^t K(t-s)* [ \int_{\mathbb R^2} |\mathbf v' | f(s) d \mathbf v' ]  ds  
 \nonumber \\
+   \tilde K(t)* \int_{\mathbb R^2} p_0 d \mathbf v'  +
 \int_0^t \tilde K(t-s)* [ \int_{\mathbb R^2}  f(s) d \mathbf v']   ds  \Big].
 \end{eqnarray}
This yields the $L^q_{\mathbf x}$ estimate on 
$\int_{\mathbb R^2} |\mathbf v| p d \mathbf v$. The proof
for  $\int_{\mathbb R^2} |\mathbf v|^2 p d \mathbf v$ is analogous.
\\


{\bf Lemma 2.7.} {\it For the solution constructed in Proposition 2.2,
let us assume that 
$|\mathbf v|^{\beta} f \in L^\infty(0,T;L^1(\mathbb R^2 \times \mathbb R^2))$
and
$|\mathbf v|^{\beta} p_0 \in L^1( \mathbb R^2 \times \mathbb R^2)$,
for $\beta \geq 0$. Then, the velocity derivatives of the solution 
$p$ of (\ref{heatlinear})-(\ref{heatdatum}) satisfy
$|\mathbf v|^{\beta} {\partial p \over \partial v_i} 
\in L^1(0,T;L^1(\mathbb R^2 \times \mathbb R^2)),$ $i=1,2$.
}

{\bf Proof.}
Differentiating identity (\ref{solint}), we find an expression for the 
derivatives of the solution:
\begin{eqnarray}
{\partial p_k \over \partial z_i}(t,\mathbf{x},\mathbf{v})= 
\int_{\mathbb{R}^2} \int_{\mathbb{R}^2} {\partial \Gamma_k \over
 \partial z_i}(t,\mathbf{x},\mathbf{v};0,\mathbf{x}',\mathbf{v}') 
 p_0(\mathbf{x}',\mathbf{v}') d \mathbf{x}' d \mathbf{v}' \nonumber \\
+ \int_0^t \int_{\mathbb{R}^2} \int_{\mathbb{R}^2} {\partial 
\Gamma_k \over \partial z_i}(t,\mathbf{x},\mathbf{v};\tau,\mathbf{x}',
\mathbf{v}') f(\tau,\mathbf{x}',\mathbf{v}')\, 
d\tau \, d \mathbf{x}' d \mathbf{v}',
\label{dsolint}
\end{eqnarray}
with $z_i=v_i$ or $z_i=x_i$.
When $\beta =0$, ${\partial p \over \partial v_i}(t)$ is an integrable
function for $t>0$ thanks to Lemma 2.1 and estimate  (\ref{bounddG}).
When $\beta >0$, we argue as in the proof of Lemma 2.6 using 
estimate (\ref{bounddG}) to obtain, for any $t \in [0,T]$:
\begin{eqnarray*}
\| |\mathbf v|^\beta {\partial p\over \partial v_i}(t)
\|_{L^1_{\mathbf x \mathbf v}} \leq 
 t^{-1/2}  C_2(T) [ \|Êp_0\|_{L^1_{\mathbf x \mathbf v}} +
 \|Ê|\mathbf v|^\betaÊp_0\|_{L^1_{\mathbf x \mathbf v}} ]
 \nonumber \\
+ C_3(T) T^{1/2} [ \|Êf\|_{L^\infty_t L^1_{\mathbf x \mathbf v}} +
 \|Ê|\mathbf v|^\betaÊf\|_{L^\infty_t L^1_{\mathbf x \mathbf v}} ].
\end{eqnarray*}

{\bf Lemma 2.8.} {\it For the solution $p$ constructed in Proposition 2.2 
when $a=a(t,\mathbf x) \in L^\infty((0,T)\times \mathbb R^2),$
$p_0 \in L^1(  \mathbb R^2 \times \mathbb R^2)$ and
$f \in L^\infty(0,T;L^1(\mathbb R^2 \times \mathbb R^2)),$
the function $\tilde p= \int_{\mathbb R^2}  p d \mathbf v$ satisfies
\begin{eqnarray} \frac{\partial}{\partial t} \tilde{p}(t,\mathbf{x}) - \sigma
 \Delta_{\mathbf{x}} \tilde{p}(t,\mathbf{x}) + 
a(t,\mathbf{x}) \tilde{p}(t,\mathbf{x}) 
= \tilde{f}(t,\mathbf{x}),  \label{integratedlinear} \\
\tilde{p}_k(0,\mathbf{x})=\tilde{p}_0(\mathbf{x}), \label{integrateddatum}
\end{eqnarray}
with source
$\tilde{f} (t,\mathbf{x})= \int_{\mathbb{R}^2} \!  
f(t,\mathbf{x},\mathbf{v}) d{\bf v}$
and initial datum
$\tilde p_0(t,\mathbf{x})=  \int_{\mathbb{R}^2} p_0
(\mathbf{x},\mathbf{v}) d\mathbf v.$}

{\bf Proof.}
 {  To justify this, we may use the integral expression (\ref{solint})
and integrate with respect to $\mathbf v$. Since the coefficients
do not depend on $\mathbf v$, the fundamental solution
is invariant by translations in $\mathbf v$ and depends on $\mathbf v - \mathbf v'$.
Thanks to (\ref{eq:integralGa})
\begin{eqnarray*}
\int_{\mathbb R^2} \Gamma(t, {\mathbf x}, {\mathbf v};  \tau, {\mathbf x}', {\mathbf v}')  
d \mathbf v  = 
\int_{\mathbb R^2} G(t - \tau, {\mathbf x} - {\mathbf x}', \mathbf{v} - \mathbf{v}')   
d \mathbf v\nonumber \\
- \int_{\tau}^t \int_{\mathbb R^2}  \left[ \int_{\mathbb R^2}
G(t-s, {\mathbf x}-{\boldsymbol \xi}, \mathbf{v}-{\boldsymbol \nu})  
 d \mathbf v \right]
a(s, {\boldsymbol \xi}) \left[
\int_{\mathbb R^2} \Gamma(s,  {\boldsymbol \xi}, {\boldsymbol \nu};   
\tau,  {\mathbf x}', \mathbf{v}') d{\boldsymbol \nu} \right]
d{\boldsymbol \xi}  ds,
\end{eqnarray*}
where $G$ is the heat kernel for diffusivity $\sigma$ in the variables
$\mathbf x, \mathbf v$. Notice that
$\int G(t - \tau, {\mathbf x} - {\mathbf x}', \mathbf{v} - \mathbf{v}')   
d \mathbf v = K(t - \tau, {\mathbf x} - {\mathbf x}')$ is the heat kernel for diffusivity $\sigma$ in the variable $\mathbf x$. Therefore,
$\tilde \Gamma = \int \Gamma(t, {\mathbf x}, {\mathbf v};  \tau, {\mathbf x}', {\mathbf v}')  d \mathbf v$ is the fundamental solution for the operator
$\frac{\partial}{\partial t} \tilde{p}(t,\mathbf{x}) - \sigma
 \Delta_{\mathbf{x}} \tilde{p}(t,\mathbf{x}) + 
a(t,\mathbf{x}) \tilde{p}(t,\mathbf{x})$. Integrating (\ref{solint})
with respect to $\mathbf v$, we conclude that $\tilde p$ is a solution
of (\ref{integratedlinear})-(\ref{integrateddatum}).

Alternatively, we may first consider smooth solutions approximating
$a$ by a mollified sequence $a_k$ (as in the proof of Proposition 2.2)
and $p_0$, $f$ by $C^\infty_c$ data with compact support $p_{k,0}$
and $f_k$. We may then
integrate equations (\ref{heatlinear})-(\ref{heatdatum}), noticing that}
$\int_{\mathbb{R}^2}  \Delta_{\mathbf v} p_k  d \mathbf v =
{\rm lim}_{R \rightarrow 0 } \int_{|\mathbf v|=R} \nabla_{\mathbf v} p_k \cdot
{\mathbf n} dS = 0,$
since the derivatives of $p_k$ are integrable functions by Lemma 2.7.
Letting $k \rightarrow \infty$, we obtain the equation for $\tilde p$
as limit of the problems for $\tilde p_k$.
 \\

{\bf Lemma 2.9.} {\it Any solution $p$ of the initial value problem 
(\ref{heatlinear})-(\ref{heatdatum}) with bounded coefficient $a\geq 0$,
initial datum $u_0 \in L^2(\mathbb{R}^2 \times \mathbb{R}^2)$
and source $f \in L^2(0,T;L^2(\mathbb{R}^2 \times \mathbb{R}^2))$ 
satisfies the energy inequality:
\begin{eqnarray}
\|p(t)\|_2^2 + 2\sigma \!\! \int_0^t \hskip -2mm 
\|Ê\nabla_{\mathbf x \mathbf v}p(s)Ê\|_2^2 ds 
\leq \|p_0\|_2^2 + 2 \! \int_0^t \hskip -2mm f(s)Êp(s) ds. 
\label{energyheat}
\end{eqnarray} 
An analogous inequality holds for solutions $\tilde p$ of 
(\ref{integratedlinear})-(\ref{integrateddatum}) and 
$c$ of (\ref{heatlinearc})-(\ref{heatdatumc}). \\
}
{\bf Proof.} 
Let us first assume that $p_0,f \in C^\infty_c$ and $a$ is replaced by a smooth mollified sequence $a_k\geq 0$ converging to $a$ in $L^\infty$ weak*. By Section 2.1,  $p_k \in C([0,T],L^2(\mathbb{R}^2 \times \mathbb{R}^2))$ and $\nabla_{\mathbf x \mathbf v} p_k \in L^2((0,T) \times \mathbb{R}^2 \times \mathbb{R}^2)$.  Using the integral expression (\ref{heatint}) and the differential equation,  $\Delta_{\mathbf x \mathbf v} p_k$ and ${\partial p_k \over \partial t}$ belong to $L^2((0,T) \times \mathbb{R}^2 \times \mathbb{R}^2)$. The equation holds in $L^2$. Multiplying the equation by $p_k$, integrating over $(0,T) \times \mathbb{R}^2 \times \mathbb{R}^2$ and integrating by parts, we find:
\begin{eqnarray}
  \int_{\mathbb{R}^2 \times \mathbb{R}^2} \hskip -8mm p_k(t)^2
  d{\mathbf x} d{\mathbf v} + 2
   \!\! \int_0^t \hskip -2mm \int_{\mathbb{R}^2 \times \mathbb{R}^2} 
 \hskip -8mm [ \sigma |\nabla_{\mathbf{x}\mathbf{v}} p_k|^2 
 + a_k |p_k|^2] d{\mathbf x} d{\mathbf v} ds = \nonumber \\
 2 \!\! \int_0^t \hskip -2mm \int_{\mathbb{R}^2 \times \mathbb{R}^2} 
 \hskip -8mm f p_k   d{\mathbf x} d{\mathbf v} ds
+   \int_{\mathbb{R}^2 \times \mathbb{R}^2} \hskip -8mm p_0^2
  d{\mathbf x} d{\mathbf v}. 
 \label{energyheatk}
\end{eqnarray}
Lemma 2.1 provides a uniform bound on $p_k$ in 
$L^2(0,T;L^2_{\mathbf x \mathbf v})$. The energy inequality 
(\ref{energyheatk}) extends this uniform bound to
$L^2(0,T;H^1_{\mathbf x \mathbf v})$. Arguing as in Proposition 2.2, 
the solutions $p_k$ of the regularized problems tend to the solution
of the problem with coefficient $a$ in $L^2(0,T;H^1_{\mathbf x \mathbf v})$ weak.
Since the limit of their norms is bounded from below by the norms
of the weak limits, taking limits in identity (\ref{energyheatk}) and neglecting
a positive term, we get inequality (\ref{energyheat}) for coefficient $a$ and smooth data of compact  support.

Now, take $C^\infty_c$ sequences $p_{k,0}$, $f_k$ converging to $p_0$,
$f$ in $L^2(\mathbb{R}^2 \times \mathbb{R}^2)$
and  $L^2(0,T;L^2(\mathbb{R}^2 \times \mathbb{R}^2))$,
respectively.  Let $p_k$ be the corresponding
solutions of problem (\ref{heatlinear})-(\ref{heatdatum}). 
Inequality (\ref{energyheat})
yields uniform $L^2(0,T;H^1_{\mathbf x \mathbf v})$ estimates, implying
weak convergence of a subsequence to a limit $p$ 
in $L^2(0,T;H^1_{\mathbf x \mathbf v})$. Taking limits
in the weak formulations for $p_k$, it follows that 
$p$ is a weak solution with data
$p_0$ and $f$. Taking limits in the energy identities for $p_k$,
we get the energy inequality (\ref{energyheat}) for $p$. \\

{\bf Lemma 2.10.} 
 {\it For the solution $p$ constructed in Proposition 2.2
 when $a=a(t,\mathbf x) \in L^\infty((0,T)\times \mathbb R^2),$
$(1+|\mathbf v|^2)p_0 \in L^1(  \mathbb R^2 \times \mathbb R^2)$
and 
$(1+|\mathbf v|^2) f \in L^\infty(0,T;L^1(\mathbb R^2 \times \mathbb R^2)),$
the function $m= \int_{\mathbb R^2} |\mathbf v|^2 p d \mathbf v$ satisfies
\begin{eqnarray} 
\frac{\partial}{\partial t} m(t,\!\mathbf{x}) \!-\! \sigma
\Delta_{\mathbf{x}} m(t,\!\mathbf{x})  +
(a(t,\!\mathbf{x})- 4 \sigma) m(t,\!\mathbf{x})   
\!=\! \hat f(t,\!\mathbf{x}),   
\label{p2integrated}
\end{eqnarray}
with integrable source
$\hat{f} (t,\mathbf{x})= \int_{\mathbb{R}^2} \! |\mathbf v|^2 
f(t,\mathbf{x},\mathbf{v}) d{\bf v}$ 
and initial datum
$\hat p_0(t,\mathbf{x})=  \int_{\mathbb{R}^2} |\mathbf v|^2p_0
(\mathbf{x},\mathbf{v}) d\mathbf v$.}

{\bf Proof.} We argue first for smooth solutions corresponding to smooth
$a$, $p_0$, $f$.
Multiplying (\ref{heatlinear}) by $|{\mathbf v}|^2$ and integrating with respect 
to ${\mathbf v}$ we obtain (\ref{p2integrated}).
Indeed, integrating by parts over balls of radius
$R$ in velocity and letting $R\rightarrow \infty$ we find:
\begin{eqnarray*} 
\sigma \int_{\mathbb R^2 \times \mathbb R^2}  \hskip -2mm
v_i^2 {\partial^2 \over \partial^2 v_i} p_k 
d\mathbf x  d\mathbf v =
- 2 \sigma \int_{\mathbb R^2 \times \mathbb R^2}  \hskip -2mm
v_i {\partial \over \partial v_i} p_k d\mathbf x  d\mathbf v =
2  \sigma \int_{\mathbb R^2 \times \mathbb R^2} \hskip -4mm
p_k d\mathbf x  d\mathbf v.
\end{eqnarray*}
The boundary integrals
$\int_{\mathbb{R}^2}  \int_{|\mathbf v|=R} 
v_i^2 {\partial \over \partial v_i} p_k  {n_i}
d\mathbf x  dS_{\mathbf v} $ and
$ \int_{\mathbb{R}^2}  \int_{|\mathbf v|=R} 
v_i  p_k  {n_i} d\mathbf x  dS_{\mathbf v} $ tend to zero as
$R$ tends to infinity as a consequence of Lemmas
2.6 and 2.7, which ensure $(1+|\mathbf v|^2) p \in 
C(0,T;L^1(\mathbb R^2 \!\times\! \mathbb R^2))$ and
 $(1+|\mathbf v|^2 ) {\partial p \over \partial v} \in 
L^1(0,T;L^1(\mathbb R^2 \!\times\! \mathbb R^2)).$
The result extends to non smooth data and coefficients by employing
approximating sequences and taking
limits, as in Proposition 2.2.
}

\section{Integrodifferential problem for the density}
\label{sec:nonlinear0}

For $k\geq 2$, we consider the iterative scheme:
\begin{eqnarray} \frac{\partial}{\partial t} p_k(t,\mathbf{x},\mathbf{v}) 
- \sigma \Delta_{\mathbf{x}\mathbf{v}} p_k(t,\mathbf{x},\mathbf{v}) + 
a_{k-1}(t,\mathbf{x}) p_k(t,\mathbf{x},\mathbf{v}) 
= f(t,\mathbf{x},\mathbf{v}),  \label{iterativelinear} \\
p_k(0,\mathbf{x},\mathbf{v})=p_0(\mathbf{x},\mathbf{v}), \label{iterativedatum}
\end{eqnarray}
for $(t, \mathbf{x},\mathbf{v}) \in [0,\infty) \times \mathbb{R}^2 \times \mathbb{R}^2$,  
with $\sigma \in \mathbb{R}^+$, $f \in L^{\infty}(0,\infty; L^{\infty}\cap L^{1} (\mathbb{R}^2 \times \mathbb{R}^2))$, $p_0 \in L^{\infty}\cap L^{1} (\mathbb{R}^2 \times \mathbb{R}^2)$, $f \geq 0$ and $p_0 \geq 0$.
The coefficient  { $a_{k-1}$ is defined as:
\begin{eqnarray}
a_{k-1}(t,\mathbf{x})= a(p_{k-1})= \int_0^t \!   ds \! \int_{\mathbb{R}^2} 
\! d{\bf v}' p_{k-1}(s,\mathbf{x},\mathbf{v}') = \int_0^t \!   ds  \, \tilde 
 p_{k-1}(s,\mathbf{x}),
\label{ak-1}
\end{eqnarray} 
for $k\geq 2$.}

We set $p_1$ equal to the solution of the heat equation obtained when $a_0=0$.  { Thanks to the integral expression (\ref{heatint}),
$p_1 \in L^{\infty}(0,T;L^{\infty}\cap L^1 (\mathbb{R}^2 \times 
\mathbb{R}^2))$ and $p_1 \in L^{\infty}(0,T;L^{\infty}_{\mathbf x}\cap  
L^{1}_{\mathbf x}(\mathbb{R}^2; L^1_{\mathbf v}(\mathbb{R}^2)))$. Additionally,  estimate (\ref{LqLrp0f}) holds.
}

An induction procedure guarantees the existence of iterates $p_k$ satisfying 
$p_k \in L^{\infty}(0,T;L^{\infty}\cap L^1 (\mathbb{R}^2 \times \mathbb{R}^2))$ 
and 
$p_k \in L^{\infty}(0,T;L^{\infty}_{\mathbf x}\cap  L^{1}_{\mathbf x}(\mathbb{R}^2; L^1_{\mathbf v}(\mathbb{R}^2)))$.   Indeed, assuming that $a_{k-1}$ is 
measurable and bounded, a unique positive solution $p_k$ exists {
in view of Proposition 2.2.}
Then, we must check that $a_k$ is a bounded function, and 
that we can construct $p_{k+1}$. { 
By Proposition 2.2 and Lemma 2.6,}
the integral expression (\ref{solint}) in terms of 
fundamental solutions satisfying (\ref{upperboundG}) 
yields the $L^1-L^{\infty}$ bounds (\ref{Linf})-(\ref{L1}) on $p_k$ 
and $L^{\infty}$ bounds on $a_k$ with constants depending on 
$\sigma$, $T$,  $\|a_{k-1}\|_{\infty}$. {
Therefore, $a_k$ is a bounded function. By induction, we can construct 
the sequence $p_k$ for all $k$ and $a_k$ is a bounded function for all 
$k$.}

Equations (\ref{iterativelinear})-(\ref{iterativedatum}) and their fundamental solutions in dimension $n=4$ ensure $L^q_{\mathbf x \mathbf v}$ regularity 
for $p_k$, { thanks to Proposition 2.2 and Lemma 2.1.} Integrating in velocity and time, we deduce $L^q_{\mathbf x}$ regularity for $\tilde p_k $ and $a_k$, { either exploting the integral equations for $p_k$ 
as in Lemma 2.6, or applying Proposition 2.2} to the differential equations (\ref{integratedlinear})-(\ref{integrateddatum}) for $\tilde p_k$ { established 
in Lemma 2.8.}

We use this iterative scheme to establish the following existence result:\\

{\bf Theorem 3.1.}
{\it There exists a nonnegative solution $p$ { of the system:}
\begin{eqnarray} 
\frac{\partial}{\partial t} p(t,\mathbf{x},\mathbf{v}) - \sigma
 \Delta_{\mathbf{x}\mathbf{v}} p(t,\mathbf{x},\mathbf{v}) + 
\int_0^t \! \! ds \!\! \int_{\mathbb{R}^2} \hskip -3mm d{\bf v}' 
p(s,\mathbf{x},\mathbf{v}') p(t,\mathbf{x},\mathbf{v}) 
= f(t,\mathbf{x},\mathbf{v}),  
\label{alinear} \\
p(0,\mathbf{x},\mathbf{v})=p_0(\mathbf{x},\mathbf{v}), 
\label{adatum}
\end{eqnarray}
satisfying
\begin{eqnarray*}
p \in L^2(0,T;H^1(\mathbb{R}^2 \times \mathbb{R}^2)) \cap
L^{\infty}(0,T;L^{\infty}\cap L^1 (\mathbb{R}^2 \times \mathbb{R}^2)),\\
p \in L^{\infty}(0,T;L^{\infty}_{\mathbf x}\cap 
L^{1}_{\mathbf x}(\mathbb{R}^2; L^1_{\mathbf v}(\mathbb{R}^2))),
\end{eqnarray*}
if 
$f \in L^{\infty} (0,T; L^{\infty} \cap L^{1} \cap H^1 (\mathbb{R}^2 \times \mathbb{R}^2))$, 
$p_0 \in L^{\infty}(0,T;L^{\infty}\cap L^1 \cap H^1 (\mathbb{R}^2 \times \mathbb{R}^2))$, 
$f\in L^{\infty} (0,T; L^{\infty}_{\mathbf x}\cap L^1_{\mathbf x}(\mathbb{R}^2; L^1_{\mathbf v}(\mathbb{R}^2))),$
$p_0 \in L^{\infty}_{\mathbf x}\cap L^1_{\mathbf x}(\mathbb{R}^2; 
L^1_{\mathbf v}(\mathbb{R}^2)),$
$ f  \geq 0, p_0 \geq 0$ and $\sigma \in \mathbb{R}^+$. This solution is 
unique and its norms are bounded in terms of the norms of the data.\\
}

We detail the steps of the proof below. After establishing a priori bounds of 
$p_k$, $a_k$, we will pass to the limit in (\ref{iterativelinear}), obtaining a solution of (\ref{alinear})-(\ref{adatum}). We will show that this solution is 
unique and study its regularity. Let us collect first the relevant a priori bounds.

\subsection{A priori bounds}
\label{sec:bounds0}

Fundamental solutions provide existence, positivity and basic regularity. However, $\|a_{k-1}\|_{\infty}$ affects the estimates in a way difficult to control. Uniform $L^q$ bounds on $p_k$ follow from comparison principles.
Since the fundamental solution of (\ref{iterativelinear}) is positive and $p_0, f \geq 0$, we have $p_k, a_{k} \geq 0$ for all $k$. { By Lemmas 2.4 and 2.5}, $p_k$ is bounded from above by the solution of the heat equation  with the same data, that is,  $p_1$. For $k\geq 2$, we have:
\begin{eqnarray}
{
0 \leq p_k(t)   \leq  p_1(t), \quad t \in [0,T].}  
\label{comparisonpkp1}
\end{eqnarray}
Since the data are integrable and bounded, $p_1$ satisfies the $L^r-L^q$ estimates (\ref{LqLrp0f}). This yields uniform bounds for $p_k$ and $a_k$. { Indeed, combining (\ref{comparisonpkp1})
and (\ref{LqLrp0f}), we get:}
\begin{eqnarray}
\|p_k (t)\|_{L^q_{\mathbf x \mathbf v}}  \leq  
\|Êp_1(t) \|_{L^q_{\mathbf x \mathbf v}} \leq 
{
C(T,\|p_0 \|_{L^q_{\mathbf x \mathbf v}},
\|f \|_{L^\infty_t L^q_{\mathbf x \mathbf v}}),}  \label{Lqk0}
\end{eqnarray}
for $t \in [0,T]$, when $k \geq 2$, $T>0$,  and $1\leq q \leq \infty$. 
{
Integrating (\ref{comparisonpkp1}) with respect to  velocity and time,
and applying Lemma 2.6 to $p_1$, we find:}
\begin{eqnarray}
0  \leq a_{k-1}(t) \leq a_1(t) \leq \|p_1\|_{L^{1} (0,t; L^{\infty}_{\mathbf x}(\mathbb{R}^2; L^1_{\mathbf v}(\mathbb{R}^2)))}
\leq { C(T,\|\tilde p_0 \|_{L^\infty_{\mathbf x}},
\|\tilde f \|_{L^\infty_t L^\infty_{\mathbf x}}),} 
\label{Linfak0} 
\end{eqnarray}
for $\mathbf{x} \in \mathbb{R}^2$, $t \in [0,T]$, when $k\geq 2$, $T>0$.

Once we have obtained uniform bounds on $p_k$ and $a_k$, the heat operator provides uniform bounds on the derivatives of $p_k$.
The starting function for the iteration $p_1$ is the solution of a problem 
for a heat equation with $L^1\cap L^{\infty}$ data. Using its integral expression in terms of heat kernels, 
\begin{eqnarray}
\nabla_{\mathbf x \mathbf v}  p_1(t)= G(t)* \nabla_{\mathbf x \mathbf v}  p_0 + \int_0^t  \nabla_{\mathbf x \mathbf v}  G(t-\tau)*f(\tau)  d\tau, 
\label{dheatint}
\end{eqnarray}
its derivatives are bounded in terms of the $L^q$
norms of the initial and source data:
\begin{eqnarray} 
\| \nabla_{\mathbf x \mathbf v}  p_1(t) \|_q 
\leq \| \nabla_{\mathbf x \mathbf v} p_0 \|_q  
+  2 t^{1/2} {\rm max}_{s\in [0,t]} \|f(s)\|_q, \quad t \in [0,T]. 
\label{dLqLrp0f}
\end{eqnarray}
We have only assumed that $p_0 \in H^1$, $H^1$ being the usual Sobolev space. Therefore, we set $q=2$.
For any $k\geq 2$, $p_k$ is a solution of a heat equation with a source term $a_{k-1}p_k$, which we have bounded in $L^{\infty}_t(L^q_{\mathbf x \mathbf v})$,  $1\leq q \leq \infty$:
\begin{eqnarray} \frac{\partial}{\partial t} p_k(t,\mathbf{x},\mathbf{v}) - \sigma
 \Delta_{\mathbf{x}\mathbf{v}} p_k(t,\mathbf{x},\mathbf{v}) = 
- a_{k-1}(t,\mathbf{x}) p_{k}(t,\mathbf{x},\mathbf{v}) 
+ f(t,\mathbf{x},\mathbf{v}). \label{iterativelinear2} 
\end{eqnarray}
This  yields $L^q_{\mathbf x \mathbf v}$ bounds on derivatives of $p_k$. Inequality  (\ref{dLqLrp0f}) with $q=2$ implies:
\begin{eqnarray} 
\| \nabla_{\mathbf x \mathbf v}  p_k (t) \|_2 &\leq& 
\| \nabla_{\mathbf x \mathbf v}  p_0 \|_2  
+  2 t^{1/2} {\rm max}_{s\in [0,t]} \left(\|a_{k-1}(s)\|_{\infty}\|p_k(s)\|_2+\|f(s)\|_2
\right)
\nonumber \\
&\leq& \| \nabla_{\mathbf x \mathbf v}  p_0 \|_2 
+  2 t^{1/2} \left( {C(p_0,f,T)}
+ \|f\|_{L^{\infty}(0,T;L^2_{{\mathbf x}{\mathbf v}})} \right)
\label{dpk}
\end{eqnarray}
in $[0,T]$, thanks to estimates { (\ref{Lqk0})
and (\ref{Linfak0}).}

As a consequence, we obtain a uniform bound on 
$\|p_k(t)\|_{L^{2}(0,T; H^1(\mathbb R^2 \times \mathbb R^2))}$ 
(which might also have been derived
from energy inequalities).  Uniform bounds on 
$\|{\partial \over \partial t} p_k(t)\|_{L^{2}
(0,T; H^{-1}(\mathbb R^2 \times \mathbb R^2))}$
follow then from (\ref{iterativelinear}). Notice that the injection 
$H^1(\Omega) \subset L^2(\Omega)$  is compact for any bounded set 
$\Omega$ \cite{brezis}. { Arguing as in the proof of Proposition 2.2,} 
the compactness results in references
\cite{lions,simon}{ (see Theorem 12.1 in \cite{lions}, Corollary 4 in \cite{simon})} 
imply the existence of a subsequence $p_{k'}$  
tending to a limit $p$ strongly on compact sets, that is, {
in ${L^{2}(0,T; L^2_{ loc}(\mathbb R^2 \times \mathbb R^2))}$,
and almost everywhere \cite{brezis}.} It also converges
weakly in  the reflexive Banach spaces in which we have uniform
bounds.  

\subsection{Convergence to a solution}
\label{sec:convergence0}

{ Thanks to the pointwise convergence obtained in the previous step we may pass to the limit in inequality (\ref{comparisonpkp1})} to obtain:
\begin{eqnarray}
0 \leq p  \leq  p_1(t) \quad \Rightarrow  \quad \|p (t)\|_q  \leq  \|Êp_1(t) \|_q
\leq C(T,\|p_0 \|_{L^q_{\mathbf x \mathbf v}},
\|f \|_{L^\infty_t L^q_{\mathbf x \mathbf v}}), \quad  \label{Lq0}
\end{eqnarray}
for $ t \in [0,T]$.
The uniform bound (\ref{comparisonpkp1}) also shows that $|p_{k'}|$ is uniformly bounded from above by a function $p_1$ belonging to $L^r(0,T;L^q_{\mathbf x \mathbf v})$ for any $1\leq r,q< \infty$. Lebesgue's dominated convergence theorem implies that $p_{k'}$ converges to $p$ in $L^r(0,T;L^q_{\mathbf x \mathbf v})$ strongly for any $1\leq r,q< \infty$.

Recall that $a_{k'-1}(t,\mathbf{x})= \int_0^t \!  ds \! \int_{\mathbb{R}^2}  d{\bf v}' p_{k'-1}(s,\mathbf{x},\mathbf{v}')$. { By inequality (\ref{comparisonpkp1}),} the integrand satisfies $0 \leq p_{k'-1} \leq p_1$. On the other hand, $p_1$ is integrable in $[0,t]\times \mathbb{R}^2 \times \mathbb{R}^2$. By Lebesgue's dominated convergence theorem, pointwise
convergence implies:
\begin{eqnarray*}
a_{k'-1}(t,\mathbf{x})= \int_0^t \!  ds \!  \int_{\mathbb{R}^2}  d{\bf v}' p_{k'-1}(s,\mathbf{x},\mathbf{v}') \longrightarrow 
a(t,\mathbf{x})= \int_0^t \!  ds \!  \int_{\mathbb{R}^2}  d{\bf v}' p(s,\mathbf{x},\mathbf{v}'),
\end{eqnarray*}
as $k$ tends to infinity, for any $t \in [0,T]$ and 
$\mathbf{x} \in \mathbb{R}^2$ fixed. 
Let us now pass to the limit in the nonlinear term $a_{k'-1} p_{k'}$. It tends to
$a p$ almost everywhere. Thanks to estimates (\ref{comparisonpkp1}) and (\ref{Linfak0}),  $0 \leq a_{k'-1} p_{k'}\leq a_1 p_1$. The upper bound $a_1p_1$
 is integrable in $[0,T]\times \mathbb{R}^2 \times \mathbb{R}^2$
because $a_1$ is bounded. Lebesgue's dominated convergence theorem yields convergence in $L^1$ and in the sense of distributions. 

As specified above, due to the uniform bounds on 
$\|p_{k'}(t)\|_{L^{2}(0,T; H^1_{{\mathbf x}{\mathbf v}})}$ 
and 
$\|{\partial \over \partial t} p_{k'}(t)\|_{L^{2}(0,T; H^{-1}_{\mathbf x \mathbf v})}$,
$p_{k'}$ tends to $p$ weakly in $L^{2}(0,T; H^1_{{\mathbf x}{\mathbf v}})$ and
${\partial \over \partial t} p_{k'}$ tends  to  ${\partial \over \partial t} p$
in $L^{2}(0,T; H^{-1}_{\mathbf x \mathbf v})$ weakly.

Let us write down the weak formulation of the
initial value problem (\ref{iterativelinear})-(\ref{iterativedatum}).
For any $\phi(t,{\mathbf x},{\mathbf v}) \in C_c^{\infty}([0,T) 
\times \mathbb{R}^2 \times \mathbb{R}^2)$,
\begin{eqnarray*}
{
\!-\!\!  \int_{\mathbb{R}^2 \times \mathbb{R}^2} \hskip -8mm 
p_0({\mathbf x},{\mathbf v})
\phi(0,{\mathbf x},{\mathbf v}) d{\mathbf x} d{\mathbf v}  
-
\!\! \int_0^T\hskip -2mm \int_{\mathbb{R}^2 \times \mathbb{R}^2}  
\hskip -8mm
f(s,\mathbf{x},\mathbf{v}) 
\phi(s,{\mathbf x},{\mathbf v}) d{\mathbf x} d{\mathbf v} ds }
\nonumber \\
= \!\! \int_0^T \hskip -2mm \int_{\mathbb{R}^2 \times \mathbb{R}^2}  
[{ {\partial \over \partial t}+} \sigma \Delta_{\mathbf{x}\mathbf{v}} \!-\!  a_{k'-1}(s,\mathbf{x})
] \phi(s,\mathbf{x},\mathbf{v}) 
p_{k'}(s,{\mathbf x},{\mathbf v}) d{\mathbf x} d{\mathbf v} ds. 
\end{eqnarray*}
Letting $k' \rightarrow \infty$ we find that $p$ is a  solution of 
(\ref{alinear})-(\ref{adatum}) in the sense of distributions and in 
$L^{2}(0,T; H^{-1}_{\mathbf x \mathbf v})$.

\subsection{Uniqueness result}
\label{sec:uniqueness0}

Let us consider first the integrated problem (\ref{integratedlinear})-(\ref{integrateddatum}) { for $\tilde p $ introduced in Lemma 2.8}:
\begin{eqnarray} \frac{\partial}{\partial t} \tilde{p}(t,\mathbf{x}) - \sigma
\Delta_{\mathbf{x}} \tilde{p}(t,\mathbf{x}) + a(t,\mathbf{x}) \tilde{p}(t,\mathbf{x}) 
= \tilde{f}(t,\mathbf{x}),  \nonumber  \\
\tilde{p}(0,\mathbf{x})=\tilde{p}_0(\mathbf{x}), \nonumber
\end{eqnarray}
where $a(t,\mathbf{x}) =\int_0^t ds \, \tilde{p}(s,\mathbf{x}) \geq 0$. 

Let us assume that we have two nonnegative solutions $\tilde{p}^{(1)}$ and $\tilde{p}^{(2)}$ belonging to  $L^{\infty}(0,T;L^{\infty}_{\mathbf x}\cap L^1_{\mathbf x})$. Set 
$\overline{p}=\tilde{p}^{(1)}-\tilde{p}^{(2)}$
and { $\overline{a}=a^{(1)}-a^{(2)}$, with 
$a^{(i)}(t,\mathbf{x}) =\int_0^t ds \, \tilde{p}^{(i)}(s,\mathbf{x}) \geq 0$,
$i=1,2$.} 
Substracting the equations for $\tilde p^{(1)}$ and $\tilde p^{(2)}$ we find:
\begin{eqnarray} \frac{\partial}{\partial t} \overline{p} - \sigma
\Delta_{\mathbf{x}} \overline{p}+ a^{(1)}\overline{p}
= -\tilde{p}^{(2)} (a^{(1)}-a^{(2)}) = -\tilde{p}^{(2)} \int_0^t (\tilde{p}^{(1)}-\tilde{p}^{(2)})(s) ds,  \label{integrateddif} \\
\overline{p}(0)=0. \label{integratedin}
\end{eqnarray}
{ Observing that $a^{(1)} \geq 0$, the energy inequality in Lemma 2.9 yields:}
\begin{eqnarray}
{1\over 2}\|\overline{p}(t)\|_2^2  + \sigma \int_0^t \! \|\nabla_{\mathbf x} \overline{p}(s)\|_2^2 ds 
\leq - \int_0^t \!\! ds \! \int_{\mathbb R^2} d\mathbf{x} \, 
\tilde{p}^{(2)}(s,\mathbf{x}) \overline{a}(s,\mathbf{x}) \overline{p}(s,\mathbf{x}).
\label{energy1}
\end{eqnarray}
%
%
Set $M_{\hat T}={\rm max}_{s\in [0,\hat T]} \|\overline{p}(s)\|_2$, $\hat T\leq T,$ and $|\tilde{p}^{(2)}(t,\mathbf{x})| \leq M$ for $0 \leq t  \leq T$ and $\mathbf{x} \in \mathbb{R}^2$. Thanks to Jensen's inequality for  convex functions { the following inequalities hold:}
\begin{eqnarray}
 | \int_0^s \overline{p}(s',{\mathbf x} ) ds' |^2  = |
{\overline a}(s,\mathbf x)|^2 \leq s \int_0^s ds' |\overline{p}(s',{\mathbf x})|^2  \nonumber \\
\Rightarrow
\| {\overline a}(s)\|_2^2 \leq s \int_0^s ds' \|\overline{p}(s')\|_2^2 \leq s^2 M_s^2.
\label{jensen1}
\end{eqnarray} 
Inserting (\ref{jensen1}) in (\ref{energy1}), we obtain:
\begin{eqnarray}
{1\over 2}\|\overline{p}(t)\|_2^2  \leq 
M (M_{\hat T})^2 {\hat T^2\over 2}, \quad t \in [0, \hat T],
\label{jensen2}
\end{eqnarray}
which implies:
\begin{eqnarray*}
(1- M \hat T^2 )M_{\hat T}^2 \leq 0.
\end{eqnarray*}
If $\hat T < {1 \over \sqrt{M}}$, this implies $M_{\hat T}=0$ and 
$\overline p= 0$ in $[0,\hat T]$. The procedure can be repeated at time $\hat T$ to get $\overline p= 0$ in $[\hat T, 2\hat T]$. Iteratively, we find $\overline p= 0$ up to time $T$, thus $a^{(1)}=a^{(2)}=a$. Then, $ {p}^{(1)}$ and $ {p}^{(2)}$ are solutions of the same linear equations, with the same initial and source data, and the same coefficient 
$a$. Therefore, they are equal { (as a consequence of either Proposition 2.2 or Lemma 2.9)} and the constructed solution is unique.

\subsection{Regularity of the solutions}
\label{sec:regularity0}

We have constructed a solution of:
\begin{eqnarray} \frac{\partial}{\partial t} p- \sigma
 \Delta_{\mathbf{x}\mathbf{v}} p=  - a \,p + f,  \quad p(0)=p_0 
 \label{diff}
\end{eqnarray}
where $a= \int_0^t \! \! ds \!\! \int_{\mathbb{R}^2} \! d{\bf v}' p(s,\mathbf{x},\mathbf{v}').$ This solution satisfies:
\begin{eqnarray*}
0 \leq p \leq p_1, \quad 0 \leq a \leq a_1,
\end{eqnarray*}
{
using first (\ref{Lq0}) and then integrating in velocity and time.}
In view of the $L^{\infty}$ bound on the coefficient $a$, 
{ the term} $a\,p$ belongs to 
$L^{\infty}_t(L^q_{\mathbf x \mathbf v})$ for any $1 \leq q \leq \infty$. 
The regularity of solutions of heat equations
implies that the derivatives of $p$ with respect to any variable
remain in $L^q_{\mathbf x \mathbf v}$. However, the $L^q_{\mathbf x \mathbf v}$ norms become singular as $t\rightarrow 0$ unless we assume regularity of the  derivatives of the initial data.

Assuming $\nabla_{\mathbf x \mathbf v} p_0 \in L^q_{\mathbf x \mathbf v}$, the integral { reformulation of the heat equation (\ref{diff}) in terms of its heat kernel $G$}
\begin{eqnarray*}
\nabla_{\mathbf x \mathbf v} p(t) = 
G(t)* \nabla_{\mathbf x \mathbf v} p_0 + \int_0^t  \nabla_{\mathbf x \mathbf v} G(t-\tau)*[f(\tau)-a(\tau)p(\tau)]  d\tau, 
\end{eqnarray*}
yields
\begin{eqnarray} 
\| \nabla_{\mathbf x \mathbf v}  p(t) \|_q \leq \| \nabla_{\mathbf x \mathbf v}  p_0 \|_q  +  2 t^{1/2} {\rm max}_{s\in [0,t]} \|f(s)- a(s) p(s)\|_q, \quad t \in [0,T]. 
\label{heatderivativenorm}
\end{eqnarray}
Otherwise, the alternative expression
\begin{eqnarray*}
\nabla_{\mathbf x \mathbf v} p(t) = 
\nabla_{\mathbf x \mathbf v} G(t)* p_0 + \int_0^t  \nabla_{\mathbf x \mathbf v} G(t-\tau)*[f(\tau)-a(\tau)p(\tau)]  d\tau, 
\end{eqnarray*}
only implies
\begin{eqnarray*} 
\| \nabla_{\mathbf x \mathbf v}  p(t) \|_q \leq t^{-1/2} \| p_0 \|_q  
+  2 t^{1/2} {\rm max}_{s\in [0,t]} \|f(s)- a(s) p(s)\|_q, \quad t \in [0,T]. 
\end{eqnarray*}
Once  estimates on the first order derivatives are available, second
order derivatives can be estimated in a similar way splitting the derivatives
between the heat kernel and the source, provided the derivatives
of $f$ also belong to $L^q_{\mathbf x \mathbf v}$, and the derivatives of $a$ 
are bounded functions. The regularity
of the time derivatives follows using the differential equation.

\section{Coupling with the diffusion equation}
\label{sec:coupling}

Let us consider now the full problem (\ref{heat1})-(\ref{heat4}) coupling the density $p$ to the { variable} $c$. The equation for the density includes now a linear source in $p$:
\begin{eqnarray*} 
\frac{\partial}{\partial t} p(t,\!\mathbf{x},\!\mathbf{v}) \!-\! \sigma
 \Delta_{\mathbf{x}\mathbf{v}} p(t,\!\mathbf{x},\!\mathbf{v}) \!+\!  
\gamma a(t,\!\mathbf{x})  p(t,\!\mathbf{x},\!\mathbf{v})  
= \alpha(c(t,\!\mathbf{x}))\rho(\mathbf{v}) p(t,\!\mathbf{x},\!\mathbf{v}).  
\end{eqnarray*}
$\rho(\mathbf{v})$ is a smooth, bounded and integrable positive function. 
This equation is coupled { with} a diffusion equation for $c$:
\begin{eqnarray*}
\frac{\partial}{\partial t}c(t,\mathbf{x})- d \Delta_{\mathbf x} c(t,\mathbf{x}) 
= - \eta c(t,\mathbf{x}) { j(t,\mathbf{x})}.
\end{eqnarray*}
Let us recall that:
\begin{eqnarray*}
\alpha(c)=\alpha_1\frac{\frac{c}{c_R}}{1+\frac{c}{c_R}}, \;
{ j(t,\mathbf{x})} \!=\! \int_{\mathbb{R}^2} \hskip -3mm 
|\mathbf{v}'| p(t,\mathbf{x},\mathbf{v}')\, d{\bf v}' , \;
a(t,\mathbf{x}) \!=\! \int_0^t  \!\! \int_{\mathbb{R}^2}  \hskip -3mm 
p(t,\mathbf{x},\mathbf{v}')\, d{\bf v}' ds.
\end{eqnarray*}
The function $c$ is expected to decay at infinity, except for a finite interval { of $x_2$} for which it tends to a constant 
{ $k_{\infty}$ as $x_1$ grows.
For any $t>0$ and $x_2 \in [a,b] \subset \mathbb R$,
\begin{eqnarray}
c(t, {x_1}, {x_2}) \rightarrow k_{\infty}  \quad \mbox{as}
 \quad {x_1} \rightarrow \infty.
\label{cinfty}
\end{eqnarray}   
That interval represents the location of a confined distant source.}
We impose the same behavior on $c(0)=c_0 \geq 0 $.
Writing $c=c_{\infty}+\hat{c}$ where $c_{\infty}$ is a solution of the heat equation with the same initial datum, $\hat{c}$ is a solution of:
\begin{eqnarray} 
\frac{\partial}{\partial t} \hat{c}(t,\mathbf{x})= 
d \Delta_{\mathbf{x}} \hat{c}(t,\mathbf{x}) 
- \eta \hat{c}(t,\mathbf{x}) { j(t,\mathbf{x})}
- \eta c_{\infty}(t,\mathbf{x}) { j(t,\mathbf{x})}, \label{eqCinf}
\end{eqnarray}
vanishing at infinity with initial datum $\hat{c}(0)=0.$ \\

The following existence and uniqueness result holds: \\

{\bf Theorem 4.1.}
{\it There exists a unique nonnegative solution $(p,c)$ of (\ref{heat1})-(\ref{heat4}) satisfying:
\begin{eqnarray*}
p \in L^2(0,T;H^1(\mathbb{R}^2 \times \mathbb{R}^2)) \cap
L^{\infty}(0,T;L^{\infty}\cap L^1(\mathbb{R}^2 \times \mathbb{R}^2)), \\
p, |{\mathbf v}|^2 p \in
L^{\infty}(0,T;L^{\infty}_{\mathbf x}\cap 
L^{1}_{\mathbf x}(\mathbb{R}^2; L^1_{\mathbf v}(\mathbb{R}^2))), \\
\hat{c}=c-c_{\infty} \in L^2(0,T;H^1(\mathbb{R}^2 )) \cap
L^{\infty}(0,T;L^{\infty}\cap L^1(\mathbb{R}^2 )), 
\end{eqnarray*}
when 
$p_0 \in L^{\infty}\cap L^1 \cap H^1(\mathbb{R}^2 \times \mathbb{R}^2)$, 
$p_0, |{\mathbf v}|^2 p_0 \in L^{\infty}_{\mathbf x}\cap L^1_{\mathbf x}(\mathbb{R}^2; L^1_{\mathbf v}(\mathbb{R}^2))$, $p_0\geq 0$ and $c_0 \in  
L^{\infty}(\mathbb{R}^2 \times \mathbb{R}^2)$, $c_0\geq 0$.    
The norms  of this solution are bounded in terms of the norms of the data.
}\\

The existence proof relies on an iterative scheme. After showing that the scheme is well defined, we obtain uniform a priori bounds on $p_k$, $a_k$, ${ j_k}$ and $c_k$. A solution of (\ref{heat1})-(\ref{heat4}) follows passing to  the limit. This solution inherits the bounds
on the iterates in terms of the norms of the data, which implies stability of the solution. Uniqueness follows from integral inequalities. We detail the proofs in the next four subsections.

%

\subsection{Iterative scheme}
\label{sec:iterativeC}

For $k\geq 2$ we consider the iterative scheme:
\begin{eqnarray} 
\frac{\partial}{\partial t} p_k(t,\mathbf{x},\mathbf{v}) - \sigma \Delta_{\mathbf{x} \mathbf{v}} p_k(t,\mathbf{x},\mathbf{v}) + \gamma a_{k-1}(t,\mathbf{x}) p_k(t,\mathbf{x},\mathbf{v})  \label{iteqp}\\ 
= \alpha(c_{k-1}(t,\mathbf{x})) \rho(\mathbf{v}) p_k(t,\mathbf{x},\mathbf{v}), \nonumber \\
\frac{\partial}{\partial t}c_{k-1}(t,\mathbf{x})= d \Delta_{\mathbf{x}}  c_{k-1}(t,\mathbf{x}) - \eta c_{k-1}(t,\mathbf{x}){ j_{k-1}(t,\mathbf{x})}. 
\label{iteqC}
\end{eqnarray}
We initialize the iteration setting $p_1=0$.  $c_1$ is the solution of (\ref{iteqC}) with initial datum $c_0$. Let us show that the iterative scheme is well defined. This follows using the fundamental solutions of the corresponding linear problems with bounded coefficients, the integral expressions for their solutions and the upper uniform bounds for the fundamental solutions involving constants depending on the $L^{\infty}$ norm of the coefficients, as we argue by induction.

Let us assume that $a_{k-1}= \int_0^t d\,s \int d{\bf v}' p_{k-1}(s,\mathbf{x},\mathbf{v}') \geq 0$ and $c_{k-1}\geq 0$ are bounded.  { By Proposition 2.2} there exists a unique positive solution $p_k$ of the initial value problem for (\ref{iteqp}), which admits an integral expression in terms of the fundamental solution {
introduced in Lemma 2.1}. This implies that $a_k $, ${ j_k}$ are bounded functions and $ a_k, j_k \geq 0$. 
Indeed, the fundamental solution is positive and bounded from above by (\ref{upperboundG}). This yields the $L^1, L^{\infty}$ bounds (\ref{Linf})-(\ref{L1}) on $p_k$ and, { by Lemma 2.6,}  $L^{\infty}$ bounds on $a_k$, ${ j_k}$ with constants depending on $\sigma$, $T$, $\gamma$, $\|a_{k-1}\|_{\infty}$, $\|\rho\|_{\infty}$ and $\alpha_1$.

For any bounded ${ j}_{k}$,  we construct a positive solution $c_{k}$ of (\ref{iteqC}) using the corresponding fundamental solution, { thanks to Lemma 2.1 and Proposition 2.3}. This yields the $L^{\infty}$ bound (\ref{Linf}) on $c_k$ with constants depending on $d$, $T$, $\eta$, and $\| { j_k} \|_{\infty}$. Moreover, $c_k$ are positive bounded functions.  

Therefore, we may repeat the procedure and construct $p_{k+1}$, $c_{k+1}$
enjoying the same properties. The iterative sequence is well defined. 

\subsection{Uniform estimates}
\label{sec:estimatesC}

Uniform estimates with respect to $k$ are a consequence of the positivity of the solutions and adequate comparison principles.

To obtain uniform $L^q$ estimates on $p_k$ we resort to the comparison principle { in Lemma 2.5.} Since $\gamma a_{k-1}p_k\geq 0$, the functions $p_k$ are bounded from above by the solution of the heat equation with the same initial datum and source 
$\alpha(c_{k-1}) \rho  \, p_{k}$:
\begin{eqnarray} 
\frac{\partial}{\partial t} u_k(t,\mathbf{x},\mathbf{v}) - \sigma
\Delta_{\mathbf{x}\mathbf{v}} u_k(t,\mathbf{x},\mathbf{v}) 
= \alpha(c_{k-1}(t,\mathbf{x})) \rho(\mathbf{v}) p_{k}(t,\mathbf{x},\mathbf{v}).  
\label{heatsource}
\end{eqnarray}
Using the integral expression (\ref{heatint}) for $u_k$ and $p_k\leq u_k$, we
get the inequality:
\begin{eqnarray}
\| p_k \|_{L^q_{\mathbf x \mathbf v}} \leq \| u_k \|_{L^q_{\mathbf x \mathbf v}}\leq 
\|p_0\|_{L^q_{\mathbf x \mathbf v}} + \alpha_1 \| \rho \|_{L^{\infty}_{\mathbf v}} \int_0^t 
\| p_{k}(s) \|_{L^q_{\mathbf x \mathbf v}} ds.
\label{heatinequality}
\end{eqnarray}
Applying Gronwall's lemma, we find:
\begin{eqnarray}
\| p_k (t) \|_{L^q_{\mathbf x \mathbf v}} \leq 
\|p_0\|_{L^q_{\mathbf x \mathbf v}} 
e^{ t \alpha_1 \| \rho \|_{\infty}}, \quad t \in [0,T],  \, 1\leq q \leq \infty.
\label{heatgronwall}
\end{eqnarray}
{ Applying again Lemmas 2.4 and 2.5,} 
$p_k \leq u_k\leq {\cal P}$, where ${\cal P}$ is a solution 
of
\begin{eqnarray} 
\frac{\partial}{\partial t} {\cal P}(t,\mathbf{x},\mathbf{v}) - \sigma
 \Delta_{\mathbf{x}\mathbf{v}} {\cal P}(t,\mathbf{x},\mathbf{v}) 
= \alpha_1 \|\rho\|_{\infty} {\cal P}(t,\mathbf{x},\mathbf{v}), 
\label{heatsource2}
\end{eqnarray}
with the same initial datum. Changing variables, it comes out that 
${\cal P} (t)= e^{\alpha_1 \|\rho\|_{\infty} t} (G(t) *_{\mathbf x \mathbf v} p_0)$, 
where $G(t)$ is the heat kernel for diffusivity $\sigma$.

{ Now, Lemma 2.8 yields} the following equation for 
$\tilde{p}_k(s,\mathbf{x})= 
\int_{\mathbb{R}^2} p_k(s,\mathbf{x},\mathbf v) d{\mathbf v}$:
\begin{eqnarray} 
\frac{\partial}{\partial t} \tilde{p}_k(t,\!\mathbf{x}) \!-\! \sigma
\Delta_{\mathbf{x}} \tilde{p}_k(t,\!\mathbf{x}) \!=\!  \alpha(c_{k-1}(t,\!\mathbf{x})) 
\!\! \int_{\mathbb{R}^2}  \hskip -3mm
d\mathbf{v} \rho(\mathbf{v}) p_k(t,\!\mathbf{v},\!\mathbf{x}) 
\!-\! { \gamma} a_{k-1}(t,\!\mathbf{x}) \tilde{p}_k(t,\!\mathbf{x})
\nonumber \\
 \leq  \alpha_1 \|\rho\|_{\infty} \tilde{p}_k(t,\!\mathbf{x}),   
\label{heatintegrated}
\end{eqnarray}
where  $a_{k-1}(t,\mathbf{x}) =\int_0^t \tilde{p}_{k-1}(s,\mathbf{x}) ds $ and 
$\tilde{p}_k(0,\mathbf{x})=\tilde{p}_0(\mathbf{x})$. 
{ By Lemma 2.4,} $\tilde{p}_k$ is bounded from above by the solution of a heat equation with source  $\alpha_1 \|\rho\|_{\infty} \tilde{p}_k(t,\mathbf{x})$.  Repeating the Gronwall argument used to estimate 
$\| p_k \|_{L^q{\mathbf x \mathbf v}}$ we find that
\begin{eqnarray}
\| \tilde{p}_k (t) \|_{L^q_{\mathbf x}} \leq 
\|\tilde p_0\|_{L^q_{\mathbf x} L^1_{\mathbf v}} 
e^{ t \alpha_1 \| \rho \|_{\infty}}, \quad t \in [0,T], \, 1\leq q \leq \infty.
\label{heatintegratedgronwall}
\end{eqnarray}
Therefore, $a_k$ is uniformly bounded in $L^{\infty}(0,T;L^q_{\mathbf x}(I\!\!R^2))$ 
for any $q\in [1,\infty]$ and any $T>0$. { Notice that estimate (\ref{heatintegratedgronwall}) would also follow directly} 
taking into account  that $p_k \leq {\cal P}$, with { ${\cal P}$ defined in (\ref{heatsource2})}. Integrating with respect to ${\mathbf v},$
we find $0 \leq \tilde{p}_k \leq \tilde{{\cal P}}$ {
and, consequently, estimate (\ref{heatintegratedgronwall}).} 
Integrating in time, we find $0 \leq a_k \leq \int_0^t \tilde{{\cal P}}(s,\mathbf{x}) ds$. 

Uniform bounds on the derivatives of $p_k$ are obtained observing that the 
functions $p_k$ solve heat equations with uniformly bounded sources in 
$L^{\infty}(0,T;L^q_{\mathbf x \mathbf v})$ for all $1 \leq q \leq \infty$:
\begin{eqnarray} 
\frac{\partial}{\partial t} p_{k}(t,\mathbf{x},\mathbf{v})  \!-\! \sigma
 \Delta_{\mathbf{x}\mathbf{v}} p_k(t,\mathbf{x},\mathbf{v})  = 
(\alpha(c_{k-1}(t,\mathbf{x}))
 \rho(\mathbf{v}) \!-\!
 { \gamma} a_{k-1}(t,\mathbf{x})) p_{k}(t,\mathbf{x},\mathbf{v}),
\label{heatsource3}
\end{eqnarray}
as in (\ref{dheatint})-(\ref{dpk}). As discussed in subsection \ref{sec:regularity0},
the derivatives of $p_k$ with respect to any variable belong to $L^q_{\mathbf x \mathbf v}$ for all $1 \leq q \leq \infty$ and all $t>0$.
When $ p_0 \in H^1(\mathbb{R}^2 \times \mathbb{R}^2)$, 
{ the partial derivatives
${\partial p_k  \over \partial z_i } \in 
L^{\infty}(0,T;L^2(\mathbb{R}^2 \times \mathbb{R}^2))$, for $z_i=x_i$ 
and $z=v_i$, $i=1,2$, thanks to inequality (\ref{heatderivativenorm}).}
A uniform bound for $p_k$ in $L^{2}(0,T,H^1_{\mathbf x \mathbf v})$ follows.
Equation (\ref{iteqp}) yields then a uniform bound on the time derivatives 
$\|{\partial \over \partial t} p_k(t)\|_{L^{2}(0,T; H^{-1}_{\mathbf x \mathbf v})}$. 

Now, we need uniform estimates on ${ j_k}$. 
{ Since we know that the $L^\infty$ norms of the coefficients 
$\alpha(c_{k-1}) \rho - \gamma a_{k-1}$ in equation (\ref{iteqp})
are uniformly bounded,
we may resort to Lemma 2.6 to obtain direct uniform estimates
on $\|{ j_k}\|_{L^\infty_t L^q_{\mathbf x}}=
\|{ j(p_k)}\|_{L^\infty_t L^q_{\mathbf x}}$ in terms
of $\|\tilde p_0\|_{L^\infty_t L^q_{\mathbf x}}$ and 
$\|j(p_0)\|_{L^\infty_t L^q_{\mathbf x}}$.}
{ Alternatively,  we can resort to uniform bounds on $m_k = \int |{\mathbf v}|^2 p_k d{\mathbf v}$ obtained from differential inequalities
provided by Lemma 2.10:} 
\begin{eqnarray} 
\frac{\partial}{\partial t} m_k(t,\!\mathbf{x}) \!-\! \sigma
\Delta_{\mathbf{x}} m_k(t,\!\mathbf{x}) \!=\! \alpha(c_{k-1}(t,\!\mathbf{x}))  
\int_{\mathbb{R}^2}  d\mathbf{v} \rho(\mathbf{v}) {
 |\mathbf v|^2p_k(t,\!\mathbf{x},\!\mathbf{v})} 
\nonumber \\
+\! ({4}\sigma \!-\! { \gamma} a_{k-1}(t,\!\mathbf{x})) 
m_k(t,\!\mathbf{x})
\leq  (\alpha_1 \|\rho\|_{\infty}+ { 4}\sigma) m_k(t,\mathbf{x}).   
\label{jintegrated}
\end{eqnarray}

In view of inequality (\ref{jintegrated}), { Lemma 2.4 and 2.5 
imply that} $m_k \geq 0$  is bounded 
from above by the solution ${\cal M}$ of a heat equation with source  
$(\alpha_1 \|\rho\|_{\infty}+ 2\sigma) m_k(t,\mathbf{x})$. 
Repeating the Gronwall argument:
\begin{eqnarray}
\| m_k (t) \|_{L^q_{\mathbf x}} \leq \||{\mathbf v}|^2 p_0\|_{L^q_{\mathbf x} 
L^1_{\mathbf v}} e^{ t (\alpha_1\|\rho\|_{\infty}+ 2\sigma)}, \quad t \in [0,T], \, 
1\leq q \leq \infty.
\label{jgronwall}
\end{eqnarray}
In fact, $0 \leq m_k \leq  {\cal M}=e^{(\alpha_1 \|\rho\|_{\infty}+2\sigma) t} 
(G(t) *_{\mathbf x} m_0)$, where $G(t)$ is the heat kernel for diffusivity 
$\sigma$.

Set ${\mathbf v}=(v_1,v_2).$ Notice that, for any $R>0$:
\begin{eqnarray} 
{ j_k} = \int_{\mathbb{R}^2} |\mathbf v|  p_k d{\mathbf v}
\leq \left[ {R} \int_{|{\mathbf v}|\leq R}  p_k d{\mathbf v} +  
\int_{|{\mathbf v}|>R} {|{\mathbf v}|^2 \over R} p_k d{\mathbf v} \right]  
\leq  ( {R} \tilde{p}_k + {1\over R} m_k ). \label{cotaj}
\end{eqnarray}
Therefore, ${ j_k}$ is uniformly bounded in 
$L^{\infty}(0,T;L^q_{\mathbf x}(\mathbb{R}^2 \times \mathbb{R}^2))$ 
for any $q\in [1,\infty]$ and $T>0$. 

Let us now obtain uniform estimates on $c_k$. The source term in 
(\ref{iteqC}) being negative, $c_{k}$ is uniformly bounded from above 
by the solution $c_{\infty}$ of the heat equation with the same initial 
datum $c_0$ { by Lemma 2.5}. Thus, $\|c_kÊ(t)\|_{\infty} \leq 
 \| c_\infty \|_{\infty} \leq \| c_0 \|_{\infty}$
for all $t \in [0,T]$. Writing down { the equation satisfied 
by $\hat{c}_k= c - c_{\infty}$,
\begin{eqnarray} 
\frac{\partial}{\partial t} \hat{c}_k(t,\mathbf{x})
- d \Delta_{\mathbf{x}} \hat{c}_k(t,\mathbf{x}) = 
- \eta  {c}_k(t,\mathbf{x}) j_k(t,\mathbf{x}), \quad
\hat{c}_k(0,\mathbf{x})=0,
\label{eqCinfk}
\end{eqnarray}
we see that} the source term is uniformly bounded in 
$L^{\infty}(0,T;L^q_{\mathbf x}(\mathbb{R}^2))$ for any $q\in [1,\infty].$ 
Then, $ \hat{c}_k, \nabla_{\mathbf x} \hat{c}_k \in 
L^{\infty}(0,T;L^2_{\mathbf x}(\mathbb{R}^2) )$ { thanks
to Lemma 2.5 and inequality (\ref{heatderivativenorm}).}
A uniform bound for $\hat{c}_k$ in $L^{2}(0,T;H^1_{\mathbf x})$ follows.
Equation (\ref{eqCinfk}) yields then a uniform estimate on the time derivatives 
$\|{\partial \over \partial t} \hat{c}_k(t)\|_{L^{2}(0,T; H^{-1}_{\mathbf x})}$. 

\subsection{Passage to the limit}
\label{sec:convergenceC}

The densities $p_k$ are uniformly bounded in $L^{2}(0,T;H^1_{\mathbf x \mathbf v})$  and their time derivatives are bounded in $L^{2}(0,T; H^{-1}_{\mathbf x \mathbf v})$. 
The { modified variables} $\hat{c}_{k}$ are bounded in 
$L^{2}(0,T;H^1_{\mathbf x})$  and their time derivatives are bounded in 
$L^{2}(0,T; H^{-1}_{\mathbf x})$.
Since the injection $H^1(\Omega) \subset L^2(\Omega)$ is compact
for any bounded $\Omega$, the classical compactness results in \cite{lions,simon} imply compactness of $p_k$  and $\hat{c}_k$ in
$L^{2}(0,T;L^2_{ loc}),$ that is, over bounded sets.
{ Arguing as in the proof of Proposition 2.2, we extract
subsequences $p_{k'}$, $\hat{c}_{k'}$ tending pointwise
and strongly to limits  $p$ and $\hat{c}$
in ${L^{2}(0,T; L^2_{{\mathbf x}{\mathbf v}}(\Omega))}$ and 
${L^{2}(0,T; L^2_{{\mathbf x}}(\omega))}$, respectively, for bounded
sets $\Omega$ and $\omega$.}
Weak convergence implies that 
$p \in L^{2}(0,T;H^1_{\mathbf x \mathbf v})$ and $\hat{c} 
\in L^{2}(0,T;H^1_{\mathbf x}).$ Weak convergences also imply that the
limits $p$ and $\hat{c} $ inherit the bounds established on the converging
sequences.

{ The subsequences $p_{k'}$ and $\hat{c}_{k'}$ converge to $p$ and $\hat{c}$ pointwise almost everywhere.} We know that $0 \leq p_{k'} \leq {\cal P}$, { where ${\cal P}$ is defined in (\ref{heatsource2})} 
and satisfies ${\cal P} \in L^r(0,T;L^q_{\mathbf x \mathbf v})$,  and any $1\leq r,q < \infty$.
Lebesgue's dominated convergence theorem implies that $p_{k'}$ converges to $p$ in $L^r(0,T;L^q_{\mathbf x \mathbf v})$. The passage to the limit in $a_{k'-1}$ and $a_{k'-1}p_{k'}$ proceeds as in Section \ref{sec:convergence0}. 

Combining continuity of $\alpha(x)$ and pointwise convergence of $c_{k'}$ to 
$\hat{c} + c_{\infty}=c$, we obtain pointwise convergence of $\alpha(c_{k'-1})$ to 
$\alpha(c).$ The positive term $\alpha(c_{k'-1}) p_{k'}$ converges almost 
everywhere to $\alpha (c) p$ and is bounded by $\alpha_1 {\cal P}$. Therefore, it 
converges strongly in any $L^r(0,T;L^q_{\mathbf x \mathbf v})$.

Recall that ${ j_{k'}}(t,\mathbf{x})= \int_{\mathbb{R}^2}  
{ |{\mathbf v} |} p_{k'}(t,\mathbf{x},\mathbf{v}) d{\mathbf v}$. 
The integrand satisfies $|{\mathbf v}| p_{k'} \leq |{\mathbf v}| {\cal P}  $,  
which is integrable over $[0,t]\times \mathbb{R}^2 \times \mathbb{R}^2$. 
By Lebesgue's dominated convergence theorem,
\begin{eqnarray*}
{ j_{k'}}(t,\mathbf{x}) \longrightarrow 
{ j}(t,\mathbf{x})=  \int_{\mathbb{R}^2} { |{\mathbf v} |}  p(t,\mathbf{x},\mathbf{v}) d{\bf v},
\end{eqnarray*}
as $k$ tends to infinity, for any $t \in [0,T]$ and $\mathbf{x} \in \mathbb{R}^2$ fixed. Let us now pass to the limit in the nonlinear term $c_{k'} {
 j_{k'}}$. It tends to $c { j}$ almost everywhere. Using (\ref{cotaj}), we
find $|c_{k'} { j_{k'}} | \leq  ( R \tilde{{\cal P}} + {1\over R} {\cal M})$, which is integrable in $[0,T]\times \mathbb{R}^2 $.
Lebesgue's dominated convergence theorem yields convergence in $L^1$ 
and in the sense of distributions. 

Passing to the limit in weak versions of the equations, we find that $(p,c)$ is a 
solution of (\ref{heat1})-(\ref{heat4}) in the sense of distributions and in 
$L^{2}(0,T; H^{-1}_{\mathbf x \mathbf v})$.

\subsection{Uniqueness}
\label{sec:uniquenessC}

{
Uniqueness follows subtracting the equations for two possible sets of solutions $p_1$, $p_2$, $c_1$, $c_2$. Let us set 
$\overline{p}=p_1-p_2$ and $\overline{c}=c_1-c_2$.  
These differences satisfy the equations:
\begin{eqnarray} 
\frac{\partial}{\partial t} \overline{p}
\!-\! \sigma \Delta_\mathbf{\mathbf x \mathbf v} \overline{p}
+ [\gamma a(p_{1}) \!-\! \alpha(c_{1}) \rho] \overline{p}  
\!=\! [-\gamma a(\overline{p}) \!+\! (\alpha(c_{1})-\alpha(c_2)) \rho] p_2,
\label{eq:pdif} \\
\frac{\partial}{\partial t} \overline{c}\!-\! d
\Delta_{\mathbf x} \overline{c} + \eta { j} (p_1) \overline{c} \!= 
- \eta  { j}(\overline{p}) c_2,
\label{eq:cdif} 
\end{eqnarray}
with $\overline{p}(0)=0$ and $\overline{c}(0)=0$. 
The mean value theorem applied to the definition $\alpha(c)$ yields:
\begin{eqnarray}
|\alpha(c_1)\!-\!\alpha(c_2)| =  \alpha_1
|{c_1 \over c_R + c_1} \!-\! {c_2 \over c_R + c_2}|
=  {\alpha_1 c_R\over (c_R+ \xi)^2} |c_1\!-\!c_2|
\leq {\alpha_1 \over c_R} |c_1\!-\!c_2|,
\label{meanalpha}
\end{eqnarray}
where $\xi \in [c_1,c_2]$. Since $c_1$ and $c_2$ are
nonnegative, $\xi  \geq 0$.

Let $\Gamma_p$ and $\Gamma_c$ denote the fundamental
solutions provided by Lemma 2.1 associated to the parabolic operators
$ \frac{\partial}{\partial t} \overline{p}
 - \sigma \Delta_\mathbf{\mathbf x \mathbf v} \overline{p}
+ [\gamma a(p_{1}) \!-\! \alpha(c_{1}) \rho] \overline{p}, $
and
$\frac{\partial}{\partial t} \overline{c}- d
\Delta_{\mathbf x} \overline{c} + \eta { j} (p_1) \overline{c},$
respectively. By the regularity properties of the solution
$(c_1,p_1)$, their coefficients are bounded functions.
Particularizing the integral expressions (\ref{solint}) provided
by Propositions 2.2 and 2.3
for $p$ and $c$, and using the upper bounds
(\ref{upperboundG}) on the fundamental solutions
$\Gamma_p$ and $\Gamma_c$, we find:
\begin{eqnarray}
\|\overline {p}(t)\|_{L^1_{\mathbf x \mathbf v}} \leq 
C_p \Big[
\gamma \|p_2\|_{L^\infty_t L^{\infty}_{\mathbf x} L^1_{\mathbf v}} 
\int_0^t \!\! ds \!\! \int_0^s \!\! d\tau 
\|\overline{p}(\tau)\|_{L^{1}_{\mathbf x \mathbf v}}
\nonumber \\
+ {\alpha_1 \|\rho\|_{\infty}\over c_R} \|p_2\|_{L^\infty_t L^{\infty}_{\mathbf x} L^1_{\mathbf v}} \int_0^t \!\! ds  \|\overline{c}(s)\|_{L^1_{\mathbf x}}
\Big], 
\label{intp1dif} \\
\|\overline {c}(t)\|_{L^1_{\mathbf x}} \leq 
C_c \eta \|c_2\|_{L^\infty_t L^{\infty}_{\mathbf x}} 
\int_0^t \!\! ds  \|j(\overline{p})(s)\|_{L^1_{\mathbf x}}, 
\label{intc1dif}
\end{eqnarray}
thanks to (\ref{meanalpha}). The constants $C_p$ and $C_c$
depend on $\sigma$, the dimension, $T$, and the $L^\infty$
norm of the coefficients.

Now, notice that $\| j(\overline{p}) \|_{L^1_{\mathbf x}} \leq 
\| |\mathbf v| \overline{p}\|_{L^1_{\mathbf x \mathbf v}}$.
This latter norm can be estimated as done in Lemma 2.6.
Using again the integral expression (\ref{solint}) for $\overline{p}$,
multiplying by $|\mathbf v|$, taking absolute values,
and integrating we obtain:
\begin{eqnarray} 
\int_{I\!\!R^2 \times I\!\!R^2}  \hskip -10mm
|\mathbf v| |\overline{p}(t, \! {\mathbf x}, \! {\mathbf v})| d \mathbf v
d \mathbf x \leq  \hskip -2mm
\int_0^t \hskip -2mm
\int_{I\!\!R^2\times I\!\!R^2\times I\!\!R^2 \times I\!\!R^2}  
\hskip - 2.2cm |\mathbf v|
\Gamma_p(t,\! {\mathbf x},\! {\mathbf v};  \!s, \!{\mathbf x}', \!{\mathbf v}')  
|f(s, \!{\mathbf x}',\! {\mathbf v}')| d{\mathbf x}' d{\mathbf v}' 
d s d \mathbf x d \mathbf v = I, \label{intvp}
\end{eqnarray}
with $f= [-\gamma a(\overline{p}) \!+\! (\alpha(c_{1})-\alpha(c_2)) \rho] p_2.$
Thanks to estimate (\ref{upperboundG}), 
\[ \Gamma_p(t, {\mathbf x}, {\mathbf v};  s, {\mathbf x}', {\mathbf v}') 
\leq C_p G(t-s,{\mathbf x}-{\mathbf x}',{\mathbf v}-{\mathbf v}'),\]
for a heat kernel $G$. Then, for $t \in [0,T]$, we have
$I\leq I_1+ I_2$, where
\begin{eqnarray}
I_1 = C_p \int_0^t \hskip -2mm (t-s)^{1\over 2} 
\hskip -2mm \int_{I\!\!R^2\times I\!\!R^2\times I\!\!R^2 \times I\!\!R^2}
\hskip -2.2cmG(t-s,{\mathbf x}-{\mathbf x}',{\mathbf v}-{\mathbf v}') 
{|\mathbf v - \mathbf v'| \over (t-s)^{1/2}}
|f(s, {\mathbf x}', {\mathbf v}')| d{\mathbf x}' d{\mathbf v}' 
d s d \mathbf x d \mathbf v,
\nonumber \\
I_2 = C_p \int_0^t \hskip -2mm
\int_{I\!\!R^2\times I\!\!R^2\times I\!\!R^2 \times I\!\!R^2}
\hskip -2.2cm
G(t-s,{\mathbf x}-{\mathbf x}',{\mathbf v}-{\mathbf v}')  
|\mathbf v'| |f(s, {\mathbf x}', {\mathbf v}')| d{\mathbf x}' d{\mathbf v}' 
d s d \mathbf x d \mathbf v. \nonumber 
\end{eqnarray}
Using the property of convolutions $\| a * b \|_1 \leq \|a \|_1 \|b\|_1$
and computing the norms involving heat kernels, we 
find:
\begin{eqnarray}
I_1 \leq  M_1 \int_0^t \| f(s) \|_{L^1_{\mathbf x \mathbf v}} ds,
\quad
I_2 \leq  M_2 \int_0^t \| |\mathbf v|f(s) \|_{L^1_{\mathbf x \mathbf v}} ds. 
\nonumber 
\end{eqnarray}
The norm $\| f \|_{L^1_{\mathbf x \mathbf v}}$ has already been estimated
in terms of the right hand side in inequality (\ref{intp1dif}) and
$\| |\mathbf v|f \|_{L^1_{\mathbf x \mathbf v}}$ is similarly bounded
taking into account the weight $|\mathbf v|$. Inserting this information in
inequality (\ref{intvp}), we get:
\begin{eqnarray}
\| |\mathbf v|\overline {p}(t)\|_{L^1_{\mathbf x \mathbf v}} \leq 
 T \Big[
\gamma [M_1\|p_2\|_{L^\infty_t L^{\infty}_{\mathbf x} L^1_{\mathbf v}} 
\!+\! M_2\||\mathbf v|p_2\|_{L^\infty_t L^{\infty}_{\mathbf x} L^1_{\mathbf v}}]
 \!\! \int_0^t \!\!\!\!  ds 
\|\overline{p}(s)\|_{L^{1}_{\mathbf x \mathbf v}}
\nonumber \\
+ {\alpha_1 \|\rho\|_{\infty}\over c_R}  [M_1\|p_2\|_{L^\infty_t L^{\infty}_{\mathbf x} L^1_{\mathbf v}} 
\!+\! M_2\||\mathbf v|p_2\|_{L^\infty_t L^{\infty}_{\mathbf x} L^1_{\mathbf v}}] \int_0^t \!\! ds  \|\overline{c}(s)\|_{L^1_{\mathbf x}}
\Big].
\label{intvp1} 
\end{eqnarray}
Notice that $ \int_0^t   ds  \int_0^s   d\tau 
\|\overline{p}(\tau)\|_{L^{1}_{\mathbf x \mathbf v}}
\leq T \int_0^t  ds 
\|\overline{p}(s)\|_{L^{1}_{\mathbf x \mathbf v}}.$

Combining inequalities (\ref{intp1dif}), (\ref{intc1dif}) and
(\ref{intvp1}), we find: 
\begin{eqnarray}
\|\overline {p}(t)\|_{L^1_{\mathbf x \mathbf v}} \leq 
A \int_0^t ds \|\overline{p}(s) \|_{L^1_{\mathbf x \mathbf v}} 
+ B  \int_0^t ds \| |\mathbf v|\overline{p}(s) \|_{L^1_{\mathbf x \mathbf v}},
\label{control1} \\
\| {|{\mathbf v}|} \overline {p}(t)\|_{L^1_{\mathbf x \mathbf v}} \leq 
\hat A \int_0^t ds \|\overline{p}(s) \|_{L^1_{\mathbf x \mathbf v}} 
+ \hat B \int_0^t ds \| |\mathbf v| \overline{p}(s) \|_{L^1_{\mathbf x \mathbf v}}.
\label{control2}
\end{eqnarray}
Setting $U(t) = {\rm max} \{\|\overline{p}(t) \|_{L^1_{\mathbf x \mathbf v}},
\| |\mathbf v| \overline{p}(t) \|_{L^1_{\mathbf x \mathbf v}}\}, 
$
we deduce from (\ref{control1})-(\ref{control2})
that $U(t)$ satisfies a Gronwall inequality of the form:
\[ U(t) \leq D \int_0^t U(s) ds\] for $t\in [0,T]$, with $D>0$. 
Therefore, Gronwall's lemma implies
$U=0$ and $p_1=p_2$ in $[0,T]$ for any $T>0$.
Then, estimate (\ref{intc1dif}) implies $c_1=c_2$.}Ê\\

{
{\bf Remark 4.2.}  The hypotheses $p_0 \in H^1_{\mathbf x \mathbf v}$
is required to establish that the iterates $p_k$ satisfy $\nabla_{\mathbf x
\mathbf v} p_k \in L^2((0,T) \times \mathbb R^2 \times \mathbb R^2)$
through integral equations for heat equations. This bound still holds
true when $p_0 \in L^2_{\mathbf x \mathbf v}$ resorting to energy
inequalities as in Lemma 2.9 instead. 
Therefore, $H^1_{\mathbf x \mathbf v}$
regularity is not needed to establish the existence and uniqueness of a solution.\\

{\bf Remark 4.3.} We have seen in the proof that the hypotheses on
$| \mathbf v |^2 p_0$ can be replaced by hypotheses on 
$| \mathbf v | p_0$  exploiting the integral equations according to
Lemma 2.6 instead of the differential inequalities provided by Lemma
2.10. Also, the uniqueness proof only needs information on
$| \mathbf v | p$.
However, arguments based on differential inequalities are more
likely to apply when trying to extend these results to bounded sets
in space $\Omega \subset \mathbb R^2$. Thus, it is worth keeping in mind
both procedures.\\

{\bf Remark 4.4.} The same existence result holds replacing $j$ by
$|\mathbf j|$, with essentially the same proof. }

\section{Discussion and future work}
\label{sec:discussion}

{
Models for angiogenesis display mathematical structures of increasing
complexity, which require the introduction of adequate strategies for
their analysis and numerical simulation. We have considered here
a simplified model, including regularizations that are frequent in
numerical approximations: replacement of Dirac measures by 
gaussians and inclusion of viscosity in degenerate directions.
We have shown that  nonnegative solutions of these regularized 
models may be constructed as limits of solutions of an iterative
scheme, obtaining stability bounds in terms of the norms of the data. 
Uniqueness conditions are also established. 
Whether regularized problems approximating measure valued 
coefficients with gaussians can be shown to effectively
converge to the original measure valued problem even in our simpler
framework is an open issue.

The main  ingredient missing in the model considered here is the transport
operator in the equation for the blood vessel density. This operator 
describes blood vessel extension in response to the chemotactic force
created by the concentration of tumor angiogenic factor. 
In principle, more realistic models including such transport operators
might be handled implementing a similar iterative procedure relying on fundamental solutions of Fokker-Planck operators for the blood vessel density, instead of fundamental solutions of 
a diffusion operator.  }

\vskip 5mm

{\bf Acknowledgements.} This work has been supported by MICINN 
and MINECO grants no. FIS2011-28838-C02-02 and No. 
MTM2014-56948-C2-1-P.

\end{document}